\documentclass{article}

\usepackage{etoolbox}
\newtoggle{long}
\togglefalse{long}

\newtoggle{ams}
\toggletrue{ams}

\newcommand{\lb}{\\ \\}

\title{Exceptional zeros and $\mathcal{L}$-invariants of Bianchi modular forms}
\date{}
\newcommand{\Addresses}{{
  \bigskip
  \footnotesize

 Daniel Barrera Salazar; Universitat Polit\`{e}cnica de Catalunya\\Campus Nord, Calle Jordi Girona, 1-3, 08034 Barcelona, Spain\par\nopagebreak
\texttt{daniel.barrera.salazar@upc.edu}

  \medskip
Chris Williams; Imperial College London\\ South Kensington Campus\\London SW7 2AZ, United Kingdom\par\nopagebreak
\texttt{christopher.williams@imperial.ac.uk}

}}

\usepackage[margin=1.45in]{geometry}
\usepackage{preamble}
\author{Daniel Barrera Salazar and Chris Williams}

\usepackage[all]{xy}

\usepackage{comment}
\usepackage{verbatim}
\iftoggle{long}{
\newenvironment{longversion}{}{} 
\newenvironment{shortversion}{\comment}{\endcomment} 

}{

\newenvironment{shortversion}{}{} 
\newenvironment{longversion}{\comment}{\endcomment} 

}

\iftoggle{ams}{
\newenvironment{amsversion}{\comment}{\endcomment} 

}{

}

\newcommand{\ind}{\mathbf{1}}

\newcommand{\tree}{\mathcal{T}}
\newcommand{\treep}{\tree_{\pri}}
\newcommand{\edges}{\mathcal{E}}

\newcommand{\edgesp}{\mathcal{E}_{\pri}}
\newcommand{\vertices}{\mathcal{V}}

\newcommand{\verticesp}{\mathcal{V}_{\pri}}
\newcommand{\Gammat}{\widetilde{\Gamma}}
\newcommand{\Omegat}{\widetilde{\Omega}}
\newcommand{\s}{\mathbb{S}}
\newcommand{\e}{e}

\newcommand{\ProjF}{\Proj(F_{\pri})}

\newcommand{\oc}{\mathrm{oc}_f}
\newcommand{\lc}{\mathrm{lc}_f}
\newcommand{\oct}{\widetilde{\mathrm{oc}}_{f,v}}
\newcommand{\lct}{\widetilde{\mathrm{lc}}_{f,\tau}}
\newcommand{\HomD}{\mathrm{Hom}(\Delta_0,V_{k,k}(\Cp)^*)}
\newcommand{\HGamma}{\mathrm{H}^1(\Gamma,\HomD)}
\newcommand{\HGammas}{\mathrm{H}^1(\Gamma,\Delta_{k,k})}
\newcommand{\rmH}{\mathrm{H}}

\newcommand{\PGLt}{\mathrm{PGL}_2}

\newcommand{\unitsize}{\#\roi_F^\times}

\newcommand{\tm}{w_{\mathrm{Tm}}}
\newcommand{\tmp}{w_{\mathrm{Tm},\pri}}

\begin{document}
\lhead{\emph{Exceptional Zeros for Bianchi Modular Forms}}
\rhead{\emph{Daniel Barrera Salazar and Chris Williams}}

\maketitle

\begin{abstract}
Let $f$ be a Bianchi modular form, that is, an automorphic form for $\GLt$ over an imaginary quadratic field $F$. In this paper, we prove an exceptional zero conjecture in the case where $f$ is new at a prime above $p$. More precisely, for each prime $\pri$ of $F$ above $p$ we prove the existence of an $\mathcal{L}$-invariant $\mathcal{L}_{\pri}$, depending only on $\pri$ and $f$, such that when the $p$-adic $L$-function of $f$ has an exceptional zero at $\pri$, its derivative can be related to the classical $L$-value multiplied by $\mathcal{L}_{\pri}$. The proof uses cohomological methods of Darmon and Orton, who proved similar results for $\mathrm{GL}_2/\Q$. When $p$ is not split and $f$ is the base-change of a classical modular form $\tilde{f}$, we relate $\mathcal{L}_{\pri}$ to the $\mathcal{L}$-invariant of $\tilde{f}$, resolving a conjecture of Trifkovi\'{c} in this case.

\end{abstract}

\section*{Introduction}
There are a wealth of results and conjectures linking important arithmetic properties of number-theoretic objects to special values of their $L$-functions. One way of attacking these conjectures is through the theory of $p$-adic $L$-functions, whose study has facilitated significant progress in examples such as the Birch and Swinnerton-Dyer conjecture. Such objects are ($p$-adic) analytic functions that interpolate special values of classical $L$-functions.
$\lb$
The interpolation property that $p$-adic $L$-functions satisfy involves a scaling factor that can, at exceptional values, vanish, forcing the $p$-adic $L$-function to vanish independently of the corresponding classical $L$-value. This happens, for example, at the value $L_p(f,0)$ if $f$ is a weight 2 newform attached to an elliptic curve $E/\Q$ at which $p$ has split multiplicative reduction. In this example, the $p$-adic $L$-function may \emph{a priori} have no meaningful relation to the central special value $L(f,1) = L(E,1)$ which plays the leading role in the Birch and Swinnerton-Dyer conjecture. In the case of newforms for $\GLt$ over $\Q$, this phenomenon was first investigated by Mazur, Tate and Teitelbaum in \cite{MTT86}, where they conjectured a relation between the derivative of the $p$-adic $L$-function and the classical $L$-function at an exceptional value. Explicitly, if $f\in S_{k+2}(\Gamma_0(N))$ is a newform with $k$ even and $p$ exactly dividing $N$, then the $p$-adic $L$-function of $f$ has an exceptional zero at the value $(\chi \tm^{k/2},k/2)$, where $ \tm$ is the Teichm\"{u}ller character and $\chi$ is any Dirichlet character with $\chi(p) = \omega$, where $-\omega$ is the Atkin-Lehner eigenvalue of $f$ at $p$. Mazur, Tate and Teitelbaum conjectured the existence of an $\mathcal{L}$-invariant $\mathcal{L}_p$ -- depending only on $f$ and $p$ -- such that 
\[\frac{d}{ds}L_p(f,\chi  \tm^{k/2},s)\big|_{s = k/2} = \mathcal{L}_p K(\chi,k/2)\Lambda(f,\chi,k/2+1)/\Omega_f^\pm,\]
where $K(\chi,k/2)$ is some explicit non-zero scalar, $\Omega_f^\pm$ is the period of $f$, and $\Lambda$ is the (normalised) $L$-function of $f$.
$\lb$
Since this conjecture was made, it has been widely studied and proved in many different ways, whilst there are a number of equivalent descriptions of the $\mathcal{L}$-invariant in terms of arithmetic data. As an example, for weight 2 it was proved by Greenberg and Stevens in \cite{GS93} and in higher weights it has been proved independently by Kato--Kurihara--Tsuji, Stevens, and Emerton using $\mathcal{L}$-invariants defined by Fontaine--Mazur, Coleman, and Breuil respectively. For a fuller description on previous work in this setting, see the introduction of \cite{BDI10}. Of particular interest for this paper are proofs of the conjecture by Darmon in \cite{Dar01} (for weight 2) and Orton in \cite{Ort04} (for higher weights) using the theory of Bruhat--Tits trees.
$\lb$
Despite the problem of exceptional zeros being well-studied for $\GLt$ over $\Q$, many of the approaches mentioned above use deep arithmetic results and structures that do not necessarily exist in the case of $\GLt$ over more general number fields, and accordingly this case has received considerably less attention. In particular, for $\GLt$ over an imaginary quadratic field, the symmetric space playing the role of the modular curve has no complex or algebraic structure, making it much harder to attach Galois representations to automorphic forms. The theory of $p$-adic variation of modular forms, which also plays a prominent role in the proof of Greenberg and Stevens, is also comparatively badly behaved; for example, there are classical Bianchi modular forms which do not vary in a $p$-adic family where the classical forms are dense.
$\lb$
Over totally real fields, the conjecture has been proved in recent years for certain weight 2 Hilbert modular forms by Mok and Spie{\ss} (see \cite{Mok09} and \cite{Spi14}). The case of higher weights has been addressed by the first author, Dimitrov and Jorza in \cite{BDJ17}. Until recently, there were no analogous results over fields with complex embeddings. The first were results of Deppe in \cite{Dep16}, where he showed that (again for certain ordinary weight 2 forms) in the cases where exceptional zeros occur, they have (at least) the expected order of vanishing. In this paper, we treat the case of arbitrary weight for $\GLt$ over an imaginary quadratic field using the methods of Darmon and Orton, which allows us to circumvent the need for the deep arithmetic structure missing in this case. We do not require an ordinarity condition, and we also make no assumption on the prime $p$. Our main result is summarised below in Theorem \ref{summarythm} of the introduction.

\subsection*{Outline of method and plan of paper}
Section 1 contains classical preliminaries. In particular, let $F$ be an imaginary quadratic field of class number 1.\footnote{This is not a serious restriction; see Section \ref{class number}.} Let $f$ be a cuspidal \emph{Bianchi modular form}, that is, a cuspidal automorphic form $f$ for $\GLt$ over $F$. Let $f$ have level $\Gamma_0(\n) \subset \SLt(\roi_F)$, with $(p)$ dividing $\n$, and weight $(k,k)$ (since cuspidal Bianchi modular forms exist only at parallel weights).\footnote{We normalise the weight in such a way that a modular elliptic curve over $F$ corresponds to a weight $(0,0)$ form. In the literature, this is sometimes referred to as `weight 2'.} Let $\Lambda$ denote the (normalised) $L$-function of $f$. A construction of $p$-adic $L$-functions for such objects is contained in \cite{Wil17} in the case where the form is an eigenform satisfying a small slope property at primes of $F$ above $p$; a summary of the relevant results of \cite{Wil17} is contained in Section 3. In this setting, the interpolation property of the $p$-adic $L$-function involves an exceptional factor for each prime $\pri$ above $p$.
$\lb$
Fix such a prime $\pri|p$ and let $\treep$ be the Bruhat--Tits tree of $\GLt(F_{\pri})$. In Section 2, we show how to vary $f$ in a `family' $\mathcal{F}$ over the edges of $\treep$. In particular, to $f$ we associate a form 
\[\mathcal{F} :\edges(\treep) \times\GLt(\A_F)\longrightarrow V_{2k+2}(\C),\]
where $V_{2k+2}(\C)$ is the space of complex homogeneous polynomials of degree $2k+2$ in two variables, in such a way that for any edge $e$, $\mathcal{F}(e,\cdot)$ is a Bianchi modular form of specified level (depending on $e$) with the same weight as $f$. We further prove (in Theorem \ref{new forms tree}) that this form is harmonic in the first variable $\treep$ if and only if $f$ is new at $\pri$. 
$\lb$
In Section 4, we define modular symbols on $\treep$ -- or, rather, families of Bianchi modular symbols indexed by $\edges(\treep)$ -- using the methods of Section 2, and define Hecke operators on the space of such symbols. The study of $\treep$ allows us to recover a classical result on the Hecke eigenvalue at $\pri$ of a $\pri$-new form as an easy corollary, and this result gives a systematic supply of exceptional zeros at central values (see Corollary \ref{exceptional zero existence}). In particular, we obtain:
\begin{mpropnum}Suppose $f$ is $\pri$-new, weight $(k,k)$ and small slope at all the primes above $p$, let $-\omega$ be the Atkin-Lehner eigenvalue of $f$ at $\pri$, and let $\varphi$ be a Hecke character for which $L(f,\varphi)$ is critical (in the sense of Deligne). Then the $p$-adic $L$-function of $f$ has an exceptional zero (for $\pri$) at $\varphi$ precisely when $k$ is even and $\varphi = \chi|\cdot|^{k/2}$, where $\chi$ is a finite order Hecke character with $\chi(\pri) = \omega$.
\end{mpropnum}
In Section 5, we utilise the overconvergent methods of \cite{Wil17} to attach a family of distributions on $\Proj(F_{\pri})$ to the form $f$. In particular, as for modular forms and modular symbols, we vary the overconvergent modular symbol attached to $f$ over the edges of $\treep$. In the case of a modular elliptic curve $E$ over $F$, this distribution has been studied before; in particular, in \cite{Tri06} Mak Trifkovic used it to construct Stark--Heegner points on $E$, whilst in the same setting the method of construction given here has been implemented independently by Xevi Guitart, Marc Masdeu and Haluk Sengun in \cite{GMS15}.
$\lb$
In Section 6, we use the distribution to define double integrals on $\uhp_{\pri}\times\uhs$ in the style of Darmon. Drawing the previous sections together, we use this to define two group cohomology classes $\lc$ and $\oc$ in $\h^1(\Gamma,\Delta_{k,k})$ (see Definition \ref{group gamma} for the definition of $\Gamma$, and Section 6 for the definition of $\Delta_{k,k})$. These classes, which depend on the choice of a complex period $\Omega_f$, are both eigenclasses for the Hecke action with the same eigenvalues as $f$. In Section 7, they are related to $L$-values using invariants attached to embeddings $F\times F \hookrightarrow \mathrm{M}_2(F)$. In particular, we prove the following:
\begin{mthmalph}Let $\chi$ be a finite order Hecke character with $\chi(\pri) = \omega$. Then there exists a map $\mathrm{Ev}_{\chi} : \h^1(\Gamma,\Delta_{k,k}) \longrightarrow \Cp$ such that
\[\mathrm{Ev}_{\chi}(\oc) = K(\chi,k/2+1)\Lambda(f,\chi, k/2+1)/\Omega_f\]
for some explicit (non-zero) factor $K(\chi,k/2+1)$, and
\[\bigg[\prod_{\mathfrak{q}}Z_{\mathfrak{q}}\bigg]\mathrm{Ev}_{\chi}(\lc) = \frac{\partial}{\partial s_{\pri}}L_p(f,\chi\tm^{k/2},\mathbf{s})\big|_{\mathbf{s} = k/2},\]
where the product is over all $\mathfrak{q}|p$ such that $\mathfrak{q} \neq \pri$, and where $Z_{\mathfrak{q}}$ is the exceptional factor at $\mathfrak{q}$.
\end{mthmalph}
In Section 8, we show that in fact, the Hecke eigenspace in which $\oc$ and $\lc$ both live is in fact one-dimensional. In particular, after proving a non-vanishing result for $\oc$, we see that there exists some element $\mathcal{L}_{\pri} \in \Cp$ such that $\lc = \mathcal{L}_{\pri}\oc$. Combined with the above, this is enough to prove the following, which is the main result of this paper.
\begin{mthmalph}\label{summarythm}
Let $f$ be as above. Let $\chi$ be a finite order Hecke character with $\chi(\pri) = \omega$, so that $L_p(f,\chi \tm^{k/2},\mathbf{k/2}) = 0$ is an exceptional zero (where $\tm$ is the Teichm\"{u}ller character). Then there exists an $\mathcal{L}$-invariant $\mathcal{L}_{\pri} \in \Cp$, depending only on $f$ and $\pri$, such that
\[\frac{\partial}{\partial s_{\pri}}L_p(f,\chi \tm^{k/2},\mathbf{s})\big|_{\mathbf{s} = k/2} = \mathcal{L}_{\pri}\cdot\bigg[\prod_{\mathfrak{q}}Z_{\mathfrak{q}}\bigg]K(\chi,k/2+1)\Lambda(f,\chi,k/2+1)/\Omega_f.\]
\end{mthmalph}
In Section 10, we discuss the case where $f$ is the base-change of a classical modular form $\tilde{f}$, and explicitly relate the $\mathcal{L}$-invariants of $f$ and $\tilde{f}$. This is enough to resolve a conjecture of Trifkovi\'{c} in the base-change case. We conclude in Section 11 by discussing possible further generalisations of these results; in particular, we discuss the case of higher class number (which has been omitted in this paper for expositional and notational reasons), possible results on the second derivative (to remove the remaining exceptional factor in the case where $p$ is split), and the conjectural form $\mathcal{L}_{\pri}$ might take in terms of arithmetic data, a question the authors hope to address in future work.

\subsection*{Comparison to relevant literature} The results of this paper fit into a body of recent work in complimentary directions. The work of Spie{\ss} (in \cite{Spi14}) and the work of the first author, Dimitrov and Jorza (in \cite{BDJ17}) on the Hilbert case has been mentioned above, as has Deppe's work (in \cite{Dep16}) on the case of weight 2 for general number fields. Exceptional zero phenomena have also been studied in the anticyclotomic setting by Molina (in \cite{Mol15}, for $\GLt$ in the totally real case) and Bergunde (in \cite{Ber17}, for non-split quaternion algebras over general number fields). Finally, the behaviour of $\mathcal{L}$-invariants under abelian base-change, again in the case of weight 2 for $\GLt$, has recently also been studied in work of Gehrmann in \cite{Geh17}.

\subsection*{Acknowledgements} The authors would like to thank Henri Darmon, who initially suggested this problem to us. The project has also benefited from helpful discussions with Tobias Berger, Kevin Buzzard, John Cremona, Xevi Guitart, David Loeffler, Marc Masdeu, and Aurel Page. The first author has received funding from the European Research Council (ERC) under the European Union's Horizon 2020 research and innovation programme (Grant Agreement No.\ 682152). The second author would like to thank Victor Rotger, whose generous financial support enabled him to visit Barcelona to complete this project. Finally, we would like to thank the anonymous referee for their comments and corrections following their careful reading of the submitted version.

%
%
\section{Preliminaries}\label{notation}

\subsection{Basic notation}\label{basic definitions}
Let $F$ be an imaginary quadratic field with ring of integers $\roi_F$, different $\mathcal{D}$ and discriminant $-D$. Let $p$ be a rational prime. If $\pri$ is a prime of $F$ above $p$, write $F_{\pri}$ for the completion of $F$ at $\pri$, write $\roi_{\pri}$ for the ring of integers in $F_{\pri}$ and fix a uniformiser $\pi$ at $\pri$. We also write $e_{\pri}$ for the ramification index of $\pri$ over $p$. In the sequel, $F$ will be taken to have class number 1, in which case we assume $\pi$ is an element of $\roi_F$ that generates $\pri$. Denote the adele ring of $F$ by $\A_F = F_\infty \times \A_F^f$, where $F_\infty$ denotes the infinite adeles and $\A_F^f$ the finite adeles. Furthermore, define $\widehat{\roi_F} \defeq \roi_F\otimes_{\Z}\widehat{\Z}$ to be the finite integral adeles. Throughout, we'll fix a prime $\pri$ above $p$ and a level $\n = \pri\m$, where $\m$ is coprime to $\pri$ and divisible by all the other primes above $p$. For an ideal $\ff \subset \roi_F$, let $\mathrm{Cl}_F(\ff)$ denote the ray class group of $F$ modulo $\ff$.
$\lb$
Let $k \geq 0$ be an integer, and for any ring $R$, let $V_{2k+2}(R)$ denote the ring of polynomials over $R$ of degree at most $2k+2$. We also define $V_0(R) = R$. Note that $V_{2k+2}(\C)$ is an irreducible complex right representation of $\SUt(\C)$, and denote the corresponding antihomomorphism by $\rho':\SUt(\C) \rightarrow \mathrm{GL}(V_{2k+2}(\C))$. Finally, define an antihomomorphism $\rho : \SUt(\C) \times \C^\times \rightarrow \mathrm{GL}(V_{2k+2}(\C))$ by $\rho(u,z) = \rho'(u)|z|^{-k}$. 
\subsection{Bianchi modular forms}\label{BMF}

A \emph{Bianchi modular form} is an automorphic form for $\GLt$ over an imaginary quadratic field. We will give only a very brief description of the theory; for a more detailed exposition of the literature, see \cite{Wil17}. First, we fix a level group.
\begin{mdef}\label{omega}
Define $\Omega_0(\n) \defeq \left\{\matrd{a}{b}{c}{d} \in \GLt(\widehat{\roi_F}): c \in \n\widehat{\roi_F}\right\}.$
\end{mdef}
\begin{mdef}We say a function $f:\GLt(\A_F) \rightarrow V_{2k+2}(\C)$ is a \emph{cusp form of weight $(k,k)$ and level $\Omega_0(\n)$} if it satisfies:
\begin{itemize}
\item[(i)]$f(zgu) = f(g)\rho(u,z)$ for $u \in \SUt(\C)$ and $z \in Z(\GLt(\C)) \cong \C^\times,$
\item[(ii)]$f$ is right-invariant under the group $\Omega_0(\n)$,
\item[(iii)] $f$ is left-invariant under $\GLt(F)$,
\item[(iv)] $f$ is an eigenfunction of the operator $\partial$ with eigenvalue $(k^2/2 + k)$, where $\partial/4$ denotes a component of the Casimir operator in the Lie algebra $\mathfrak{s}\mathfrak{l}_2(\C)\otimes_{\R}\C$ (see \cite{Hid93}, section 1.3), and 
\item[(v)] $f$ satisfies the cuspidal condition that for all $g \in \GLt(\A_F)$, we have 
\[\int_{F\backslash\A_F}f(ug)du = 0,\]
where we consider $\A_F$ to be embedded inside $\GLt(\A_F)$ by the map sending $u$ to $\smallmatrd{1}{u}{0}{1}$, and $du$ is the Lebesgue measure on $\A_F$.
\end{itemize}
Denote the space of such functions by $S_{k,k}(\Omega_0(\n))$.
\end{mdef}
\begin{longversion}
There is a good theory of Hecke operators on spaces of Bianchi modular forms, and there are multiplicity one theorems for eigenspaces. In particular, the Atkin-Lehner theory of oldforms and newforms passes over to this setting in a natural way; we will make this more precise in the sequel. 
\end{longversion}
\begin{mrem}
Note that if $f \in S_{k,k}(\Omega_0(\n))$, then $f(gh) = f(g)$ for all $h\in Z(\GLt(\A_F^f)).$ Hence we can consider $f$ as a function on $\PGLt^f(\A_F) \defeq \GLt(F_\infty) \times \PGLt(\A_F^f)$.
\end{mrem}

\subsection{The $L$-function of a Bianchi modular form}
\begin{longversion}
We recap the definitions of Hecke operators, via double coset operators. Let $I \subset \roi_F$ be an ideal, recall the definition of $x_I$ from Section \ref{basic definitions}, and choose a (finite) set of representatives $\delta_j^{I}$ such that
\[\Omega_0(\n)\matrd{1}{0}{0}{x_I}\Omega_0(\n) = \coprod_{j}\Omega_0(\n)\delta_{j}^{I}.\]
\begin{mdef}\label{Hecke}
The Hecke operator at $I$ is defined by
\[f|T_{I}(g) = \sum_j f\left(g\delta_j^{I}\right).\]
When $I$ is not coprime to $\n$, we instead write $U_{I}$ for this operator.
\end{mdef}
\end{longversion}
\begin{shortversion}There is a good theory of Hecke operators (indexed by ideals of $\roi_F)$ on Bianchi modular forms.
\end{shortversion} Let $f$ be a cuspidal Bianchi modular form that is an eigenform for all of the Hecke operators, and for any non-zero ideal $I \subset \roi_F$, write $f|T_{I} = \lambda_{I}f.$ 
\begin{mdef}
The (normalised) \emph{$L$-function of $f$} is the function
\[
\Lambda(f,\varphi) \defeq \frac{\Gamma(q+1)\Gamma(r+1)}{(2\pi i)^{q+r+2}}\sum_{I \subset \roi_F, I\neq 0} \lambda_I\varphi(I) N(I)^{-1},
\]
where $\varphi$ is a Hecke character of infinity type $(q,r)$ and $\Gamma$ denotes the usual $\Gamma$-function.
\end{mdef}
The `critical' values of this $L$-function can be controlled; in particular, we have the following:
\begin{mprop}\label{period}
There exists a period $\Omega_{f} \in \C^\times$ and a number field $E$ such that, if $\varphi$ is a Hecke character of infinity type $0 \leq (q,r) \leq (k,k)$, with $q,r \in \Z$, we have
\[\frac{\Lambda(f,\varphi)}{\Omega_{f}} \in E(\varphi),\]
where $E(\varphi)$ is the number field over $E$ generated by the values of $\varphi$.
\end{mprop}
\begin{proof}
See \cite{Hid94}, Theorem 8.1.
\end{proof}
In the case of newforms, as treated in the sequel, exceptional zeros occur at very specific Hecke characters of the form $\chi|\cdot|^r$, where $\chi$ is a finite order Hecke character and $|\cdot|$ is the adelic norm. In this case, we will write 
\[\Lambda(f,\chi,r+1) \defeq \Lambda(f,\chi|\cdot|^r)\]
for convenience and to ensure consistent notation in later sections.


%
%
\section{Forms on the Bruhat--Tits tree}
In this section, we develop the theory of Bianchi modular forms on the Bruhat--Tits tree. We will give the theory in quite general terms, since it is not significantly more difficult to describe the results in adelic language instead of classical language, and the more general results may be of independent interest. In later sections, we will specialise to the case where the class number is 1.

\subsection{The Bruhat--Tits tree}
\begin{mdef}
Let $\pri$ be a prime of $F$ above $p$. We denote by $\treep$ the Bruhat--Tits tree for $\GLt(F_{\pri})$, that is, the connected tree whose vertices are homothety classes of lattices in a two-dimensional $F_{\pri}$-vector space $V$. Two vertices $v$ and $v'$ are joined by an edge if and only if there are representatives $L$ and $L'$ respectively such that
\[\pi_{\pri}L' \subset L \subset L'.\]
Each edge comes equipped with an orientation. Denote the set of (oriented) edges of $\treep$ by $\mathcal{E}(\treep)$ and the set of vertices by $\mathcal{V}(\treep).$
\end{mdef}
Define the \emph{standard vertex} $v_*$ to be the vertex corresponding to the lattice $\roi_{\pri}\oplus\roi_{\pri}$, and the \emph{standard edge} $e_*$ to be the edge joining $v_*$ to the vertex corresponding to $\roi_{\pri}\oplus\pi_{\pri}\roi_{\pri}$. There is a natural notion of \emph{distance} between two vertices $v$ and $v'$, and we say a vertex $v$ is \emph{even} (respectively \emph{odd}) if the distance between $v$ and $v_*$ is even (respectively odd). Each (oriented) edge has a source and a target vertex, and we say such an edge is even if its source is. Write $\vertices^+(\treep)$ and $\edges^+(\treep)$ (respectively $\vertices^-(\treep)$ and $\edges^-(\treep)$) for the set of even (respectively odd) vertices and edges respectively.
$\lb$
There is a natural transitive action of $\PGLt(F_{\pri})$ on $\treep$. We can extend this action to a larger group.
\begin{mdef}\label{group omega}
Let $\m$ be an ideal of $\roi_F$ coprime to $\pri$.
\begin{itemize}
\item[(i)] For $v$ a finite place of $F$, define 
\[R_0(\m)_v \defeq \left\{\matr\in\mathrm{M}_2(\roi_v): c \equiv 0 \newmod{\m}\right\}.\]
\item[(ii)] Let $R = R_0(\m) \defeq \left\{\gamma \in \mathrm{M}_2\left(\A_F^f\right): \gamma_{v} \in R_0(\m)_{v}\text{ for }v\neq \pri, \gamma_{\pri}\in\mathrm{M}_2(F_{\pri}) \right\}.$
\item[(iii)]Let $\Omegat$ denote the image of $R^\times$ in $\mathrm{PGL}_2\left(\A_F^f\right)$.
\item[(iv)] Finally, let $\Omega \subset \widetilde{\Omega}$ denote the intersection $\Omega \defeq \mathrm{PGL}_2^+\left(\A_F^f\right)\cap\Omegat,$ where 
\[\PGLt^+\left(\A_F^f\right) \defeq \left\{\gamma \in \PGLt\left(\A_F^f\right):v_{\pri}(\det(\gamma_{\pri})) \equiv 0\newmod{2} \text{ for all }\pri|p\right\}.\]
\end{itemize}
\end{mdef}
These groups act on $\treep$ via projection to $\mathrm{PGL}_2(F_{\pri})$. 
\begin{mprop}\begin{itemize}
\item[(i)]The group $\Omegat$ acts transitively on the sets $\vertices(\treep)$ and $\edges(\treep)$. 
\item[(ii)]
The group $\Omega$ acts transitively on the sets $\edges^\pm(\treep)$ and $\mathcal{V}^\pm(\treep)$.
\end{itemize}
\end{mprop}
\begin{proof}
See \cite{Ser80}, Theorem 2 of Chapter II.1.4.
\end{proof}

A key reason for introducing the Bruhat--Tits tree is the conceptually nice description it gives of open sets in the projective line. We have a (left) action of $\PGLt(F_{\pri})$ on $\Proj(F_{\pri})$ by
\[\matr \cdot x = \frac{b+dx}{a+cx},\]
and then to $e = \gamma e_*\in\edges(\treep)$, we attach the open set 
\[
U(e) \defeq \gamma^{-1}(\roi_\pri) = \left\{x\in\Proj(F_{\pri}): \gamma x \in \roi_{\pri}\right\} \subset \Proj(F_{\pri}).
\]
\begin{longversion}
\begin{mrem}
The authors apologise for the regrettable nature of this action. In this instance it would be much nicer to consider, as in \cite{Dar01} and \cite{Ort04}, the natural action by Mobius transformations. However, much the theory of modular symbols -- including the theory over imaginary quadratic fields in \cite{Wil17} -- has been developed using this `twisted' action. Since we are unable to maintain notation that is compatible with both approaches, we have opted for this twisted action to avoid conflict when we later consider (classical and overconvergent) modular symbols on the tree. The reader should also note that our convention that $U(e) = \gamma^{-1}\roi_{\pri}$ is inverse to that of Darmon, who defines it to be $\gamma\roi_{\pri}$. Again, this difference is down to the change of action. 
\end{mrem}
\end{longversion}
\begin{shortversion}
\begin{mrem}
The authors apologise for the regrettable nature of this action, which has been chosen to be consistent with \cite{Wil17}, whilst unavoidably differing from the more natural definition of \cite{Dar01} and \cite{Ort04}. (In the sequel, $\GLt(F)$ will act on the cusps $\Proj(F)$ in this more natural way).
\end{mrem}
\end{shortversion}
\begin{mprop}\label{basis of opens}
As $e$ ranges over all edges in $\edges(\treep)$, the sets $U(e)$ form a basis of open sets in $\Proj(F_{\pri})$.
\end{mprop} 
In the sequel, we'll use the theory of `modular symbols on the tree' and this fact to construct distributions on $\Proj(F_{\pri})$ that allow us to define analogues of Darmon's double integrals in this setting.

\subsection{Bianchi forms on the tree}
Modular forms on the Bruhat--Tits tree were first defined (for classical weight 2 modular forms) in Chapter 1.1 of Darmon's seminal paper \cite{Dar01}, and the concepts are thoroughly motivated and explained in his account. The reader unfamiliar with these concepts is strongly urged to read Darmon's account, which is beautiful and well-written. In this section, we give natural analogues of his results in the imaginary quadratic setting. \\
\\
Let $R$, $\Omegat$ and $\Omega$ be as above.
\begin{mdef}
We say a function $f: \edges(\treep) \rightarrow V_{2k+2}(\C)$ is \emph{harmonic} if
\[f(\overline{e}) = -f(e) \hspace{12pt}\forall e\in\edges(\treep)\]
and if for all vertices $v\in\vertices(\treep)$, we have
\[\sum_{e: s(e) = v}f(e) = \sum_{e: t(e) = v}f(e) = 0.\]
\end{mdef}
Recall the definition $\PGLt^f(\A_F) \defeq \GLt(\C)\times\PGLt(\A_F^f)$, and note that $\Omega$ acts on $\PGLt^f(\A_F)$ by right multiplication. We define Bianchi cusp forms on $\treep$ as follows. 
\begin{mdef}
A \emph{cusp form} on $\treep\times\PGLt^f(\A_F)$ of weight $(k,k)$ for $\Omega$ is a function
\[\mathcal{F}:\edges(\treep)\times\PGLt^f(\A_F) \longrightarrow V_{2k+2}(\C),\]
such that:
\begin{itemize}
\item[(i)] $\mathcal{F}(\gamma e,g\gamma) = \mathcal{F}(e,g)$ for all $\gamma \in \Omega$.
\item[(ii)] $\f$ is harmonic as a function on $\edges(\treep)$.
\item[(iii)] For each edge $e \in \edges(\treep)$, the function
\[\mathcal{F}_{e}(g) \defeq \mathcal{F}(e,g)\]
is a cusp form of weight $(k,k)$ and level $\Omega_{e} \defeq \mathrm{Stab}_\Omega(e).$
\end{itemize}
Denote the space of such forms by $\s_{k,k}(\Omega,\treep)$.
\end{mdef}

It is natural to ask: for which $f \in S_{k,k}(\Omega_0(\pri\m))$ does there exist some $\f \in \s_{k,k}(\Omega,\treep)$ such that $\f_{e_*} = f$? The answer is:
\begin{mthm}\label{new forms tree}
The association $\f \mapsto \f_{e_*}$ defines an isomorphism
\[\s_{k,k}(\Omega,\treep) \cong S_{k,k}(\Omega_0(\pri\m))^{\pri-\mathrm{new}}.\]
\end{mthm}
We prove this in the next section. In order to do so, we introduce the following auxiliary spaces which are easier to work with, since their definition omits the harmonicity condition.
\begin{mdef}
A \emph{cusp form} on $\edges(\treep)\times\PGLt^f(\A_F)$ of weight $(k,k)$ for $\Omegat$ is a function
\[\f:\edges(\treep)\times\PGLt^f(\A_F) \longrightarrow V_{2k+2}(\C)\]
such that:
\begin{itemize}
\item[(i)] $\f(\gamma e,g \gamma) = \mathcal{F}(e,g)$ for all $\gamma \in \Omegat$.
\item[(ii)] For each $e$, the function $\f_{e}(g) \defeq \f(e,g)$ is a cusp form of weight $(k,k)$ and level $\widetilde{\Omega}_{e} \defeq \mathrm{Stab}_{\Omegat}(e).$
\end{itemize}
Denote the space of such forms by $\s_{k,k}(\Omegat,\edgesp).$
\end{mdef}
Let $f \in S_{k,k}(\Omega_0(\pri\m))$ be a Bianchi newform, and note that $f$ gives rise to a form $\f$ on $\mathcal{E}(\treep)\times\PGLt^f(\A_F)$ in a natural way. Indeed, since $\Omegat$ acts transitively on $\edges(\treep)$, any form on $\mathcal{E}(\treep)\times\PGLt^f(\A_F)$ is uniquely determined by its restriction to $\{e\}\times\PGLt^f(\A_F)$, where $e$ is any element of $\edges(\treep)$.\\
\\
An explicit check shows that $\mathrm{Stab}_{\Omegat}(e_*) = \Omega_0(\pri\m).$ Hence, writing $\f(e_*,g) = f(z,t)$ and extending $\f$ to all of $\mathcal{E}(\treep)$ using property (i), we get a well-defined map
\[S_{k,k}(\Omega_0(\pri\m)) \longrightarrow \s_{k,k}(\Omegat,\edgesp).\]
This map has a natural inverse given by $\f \longmapsto \f_{e_*}.$ We see that we have proved:
\begin{mprop}\label{iso edges}
There is a natural isomorphism $\s_{k,k}(\widetilde{\Omega},\edgesp)\cong S_{k,k}(\Omega_0(\pri\m))$ given by $\f \mapsto \f_{e_*}.$
\end{mprop}
We also need the notion of degeneracy maps at a prime $\pri$. To define these, we use the following spaces.
\begin{mdef}
A \emph{cusp form} on $\vertices(\treep)\times\PGLt^f(\A_F)$ of weight $(k,k)$ for $\Omegat$ is a function
\[\f:\vertices(\treep)\times\PGLt^f(\A_F) \longrightarrow V_{2k+2}(\C)\]
such that:
\begin{itemize}
\item[(i)] $\f(\gamma v,g \gamma) = \f(v,g)$ for all $\gamma \in \Omegat$.
\item[(ii)] For each $v\in\vertices(\treep)$, the function $\f_{v}(g) \defeq \f(v,g)$ is a cusp form of weight $(k,k)$ and level $\widetilde{\Omega}_{v} \defeq \mathrm{Stab}_{\Omegat}(v).$
\end{itemize}
Denote the space of such forms by $\s_{k,k}(\Omegat,\verticesp).$
\end{mdef}

Using the same ideas as above, replacing $e_*$ with $v_*$ and noting that $\mathrm{Stab}_{\Omegat}(v_*) = \Omega_0(\m),$ we get:
\begin{mprop}\label{iso vertices}
There is a natural isomorphism $\s_{k,k}(\Omegat,\vertices_{\pri}) \cong S_{k,k}(\Omega_0(\m))$ given by $\f \mapsto \f_{v_*}.$
\end{mprop}
\begin{longversion}
\begin{mrem}
The lack of a harmonicity condition, which we have dropped in the definitions above, makes the proof of Propositions \ref{iso edges} and \ref{iso vertices} very straightforward. This harmonicity condition will be crucial in the sequel, however, since we will use forms on the tree to define systems of distributions, and these will only be well-defined when the form is harmonic. Hence we really do need the more difficult Theorem \ref{new forms tree} to proceed.
\end{mrem}
\end{longversion}

\subsection{Proof of Theorem \ref{new forms tree}}
We now turn to the proof of Theorem \ref{new forms tree}, which we will prove by setting up two compatible exact sequences; one at the level of forms on the tree, and one at the level of classical forms. First, we focus on the tree.
\begin{mdef}\label{alpha}Let $\alpha \in \mathrm{PGL}_2(F_{\pri})\backslash\mathrm{PSL}_2(F_{\pri})$ be an element of the normaliser of $R_0(\pri)_{\pri}$, so that $\alpha e_* = \overline{e_*},$ and note that $\alpha$ switches the parity of the edges and vertices of $\treep$. 
\end{mdef}
\begin{mdef}Define a homomorphism $i : \s_{k,k}(\Omega,\treep) \longrightarrow \s_{k,k}(\Omegat,\edgesp)$ by
\[i(\f)(e,g) = \left\{\begin{array}{ll}\f(e,g) &: e\text{ even},\\
\f(\alpha e,g\alpha) &: e\text{ odd}.
\end{array}
\right.
\]
\end{mdef}
Note that by definition we've extended the $\Omega$-invariance to $\Omegat$-invariance and that $i$ is injective by the harmonicity condition.
\begin{mdef}
Define two \emph{degeneracy maps} $\pi_s,\pi_t: \s_{k,k}(\Omegat,\edgesp)\longrightarrow \s_{k,k}(\Omegat,\verticesp)$ by setting
\begin{longversion}\begin{align*}\pi_s(\f)(v,g) &\defeq \sum_{e \in \edges(\treep): s(e) = v}\f(e,g),\\
\pi_t(\f)(v,g) &\defeq \sum_{e \in \edges(\treep): t(e) = v}\f(e,g).
\end{align*}\end{longversion}
\begin{shortversion}
\[\pi_s(\f)(v,g) \defeq \sum_{e \in \edges(\treep): s(e) = v}\f(e,g)\]
(and similarly for $\pi_t$ with $t(e)$ replacing $s(e)$).
\end{shortversion}

\end{mdef}

The kernel of $\pi_s\oplus\pi_t$ is precisely the image of $\s_{k,k}(\Omega,\treep)$ under the map $i$. To see this, note that $\mathrm{Im}(i) \subset \mathrm{Ker}(\pi_s\oplus\pi_t)$ from the definition of harmonicity; conversely, suppose $\tilde{f}$ is an element of the kernel. We can directly construct $f \in \s_{k,k}(\Omega,\treep)$ with $i(f) = \tilde{f}$ by defining $f(e,g) = \tilde{f}(e,g)$ for $e$ even and extending to odd edges using harmonicity. We see that we get an exact sequence
\[0 \longrightarrow \s_{k,k}(\Omega,\treep) \longrightarrow \s_{k,k}(\Omegat,\edgesp) \labelrightarrow{\pi_s\oplus\pi_t} \s_{k,k}(\Omegat, \verticesp)\oplus \s_{k,k}(\Omegat,\verticesp).
\]

There are classical analogues of the degeneracy maps above. Write $\Omega_0(\m) = \coprod_{i\in I_{\pri}}\Omega_0(\pri\m)\gamma_i$ for a system of coset representatives $\gamma_i \in \Omega_0(\m)$. Then define:
\begin{mdef}
Let $\varphi_s$ and $\varphi_t$ denote the degeneracy maps $\varphi_s,\varphi_t: S_{k,k}(\Omega_0(\pri\m)) \longrightarrow S_{k,k}\left(\Omega_0(\m)\right)$ defined by 
\begin{align*}\varphi_s(f)(g) &\defeq \sum_{i\in I_{\pri}}f(g\gamma_i),\\
\varphi_t(f)(g) &\defeq \sum_{i\in I_{\pri}}f(g\gamma_i\alpha).
\end{align*}
\end{mdef}
\begin{mdef}
The subspace of $S_{k,k}(\Omega_0(\pri\m))$ of \emph{$\pri$-new forms} is the kernel of the map
\[S_{k,k}(\Omega_0(\pri\m)) \labelrightarrow{\varphi_s\oplus\varphi_t} S_{k,k}(\Omega_0(\m))\oplus S_{k,k}(\Omega_0(\m)).\]
\end{mdef}
\subsubsection{Relating the exact sequences}
\begin{longversion}The key step of the proof of Theorem \ref{new forms tree} is the following.\end{longversion}
\begin{mlem}
The following diagram of exact sequences commutes:
\[\begin{diagram}
 0   &\rTo&   \s_{k,k}(\Omega,\treep) &\rTo^{i}&  \s_{k,k}(\Omegat,\edgesp) &&\rTo^{\pi_s\oplus\pi_t}&& \s_{k,k}(\Omegat, \verticesp)^2 \\
&&\dTo&&\dTo&&&&\dTo\\
 0   &\rTo&   S_{k,k}(\Omega_0(\pri\m))^{p-\mathrm{new}} &\rTo&  S_{k,k}(\Omega_0(\pri\m)) &&\rTo^{\varphi_s\oplus\varphi_t}&&  S_{k,k}(\Omega_0(\m))^2 \\
\end{diagram},
\]
where the vertical arrows are the maps of the form $\f \mapsto \f_{e_*}$ (for the first two) and $\f \mapsto \f_{v_*}$ (for the third).
\end{mlem}
\begin{proof}
We've already shown that the sequences are exact. The first square commutes by the definition of $i$, since $e_*$ is even, and hence $[i(\f)]_{e_*} = \f_{e_*}$. The second square commutes since the sets
\[\{\gamma_i e_* : i \in I_{\pri}\} \hspace{6pt}\text{and}\hspace{6pt} \{\gamma_i \alpha e_* : i \in I_{\pri}\}\]
form a complete set of edges with source and target $v_*$ in $\treep$.
\end{proof}
Theorem \ref{new forms tree} now follows from the 5-lemma applied to this diagram in conjunction with Propositions \ref{iso edges} and \ref{iso vertices}.

%
%

\section{The $p$-adic $L$-function of a Bianchi modular form}
We briefly recap the results of \cite{Wil17}, which will be used heavily in the sequel. For the rest of the paper, we will make the following simplifying assumption (which was \emph{not} made in \cite{Wil17}):
\begin{mass}
The imaginary quadratic field $F$ has class number 1.
\end{mass}
\subsection{Bianchi modular symbols}\label{BMF}
Define
\[\Gamma_0(\n) \defeq \SLt(F) \cap [\GLt(\C)\Omega_0(\n)].\]
\begin{mdef}
Let $\Delta_0 \defeq \mathrm{Div}^0(\Proj(F))$ denote the space of `paths between cusps' in $\uhs$, and let $V$ be any right $\SLt(F)$-module. For a subgroup $\Gamma\subset\SLt(F)$, denote the space of \emph{$V$-valued modular symbols for $\Gamma$} to be the space
\[\symb_{\Gamma}(V) \defeq \mathrm{Hom}_\Gamma(\Delta_0,V)\]
of functions satisfying the $\Gamma$-invariance property that
\[(\phi|\gamma)(D)\defeq \phi(\gamma D)|\gamma = \phi(D)\hspace{12pt} \forall D\in\Delta_0, \gamma\in\Gamma,\]
where $\Gamma$ acts on the cusps by $\smallmatrd{a}{b}{c}{d}\cdot r = (ar+b)/(cr+d).$ For $r,s \in \Proj(F)$, to ease notation we will henceforth write $\phi\{r-s\}$ for $\phi(\{r\}-\{s\})$.
\end{mdef}
\begin{mdef}For a ring $R$, let $V_{k,k}(R) \defeq V_k(R)\otimes_{R}V_k(R)$. (We think of $V_{k,k}$ as polynomials on $\roi_F\otimes_{\Z}\Zp$ that have degree at most $k$ in each variable).
\end{mdef}
This space has a natural left action of $\GLt(R)^2$ induced by the action of $\GLt(R)$ on each factor by
\[\matr\cdot P(z) = \frac{(a+cz)^k}{(ad-bc)^{k/2}} P\left(\frac{b+dz}{a+cz}\right),\]
which is well-defined since we took $k$ to be even. This induces a right action on the dual space $V_{k,k}(R)^* \defeq \mathrm{Hom}(V_{k,k}(R),R).$ 
\begin{mrems}
\begin{itemize}
\item[(i)]Note the factor of the determinant, which was not needed in \cite{Wil17}; this ensures that the centre of $\GLt(R)$ acts trivially. This difference means that the Hecke operators we consider in this paper are scalar multiples of those in \cite{Wil17}.
\item[(ii)]We can see any subgroup of $\GLt(F)$ as acting on $V_{k,k}(\C)$ via the natural embedding of $F\hookrightarrow \C\times\C$, where the first factor differs from the second by complex conjugation. In particular, we obtain an action of $\Gamma_0(\n)$ on $V_{k,k}(\C)$.
\end{itemize}
\end{mrems}
\begin{longversion}
This also gives a natural Hecke action on the space $\symb_{\Gamma_0(\n)}(V_{k,k}(\C)^*)$, defined essentially as in Definition \ref{Hecke}. Now and for the rest of the paper, we fix more concrete choices of $\delta_{j}^{I}$.
\end{longversion}
\begin{mdef}\label{delta}
Choose $\beta_I$ to be a generator of the ideal $I$, then choose $\delta_{j}^I \in \GLt(F)\cap \mathrm{M}_2(\roi_F)$ such that
\[\Gamma_0(\n)\matrd{1}{0}{0}{\beta_I}\Gamma_0(\n) = \coprod_j\Gamma_0(\n)\delta_j^{I}.\]
Then define the Hecke operator $T_I$ by
\[(\phi|T_I)\{r-s\}(P) \defeq N(I)^{k/2}\sum_j \phi\{\delta_j^I r - \delta_j^I s\}(\delta_j^I\cdot P).\]
We write $U_I$ if $I$ is not coprime to $\n$. In the case $I = \pri = (\pi)$, we can choose the representatives to be $\smallmatrd{1}{a}{0}{\pi}$ as $a$ ranges over classes $\newmod{\pri}$.
\end{mdef}
\begin{mprop}\label{attach modular symbol}Let $f \in S_{k,k}(\Omega_0(\n))$ be a cuspidal Bianchi modular form.
\begin{itemize}
\item[(i)] To $f$, one can associate a canonical element $\widetilde{\phi}_{f} \in \symb_{\Gamma_0(\n)}(V_{k,k}(\C)^*)$ in a way that is equivariant with respect to the action of Hecke operators.
\item[(ii)] The symbol $\phi_f \defeq \widetilde{\phi}_{f}/\Omega_f$, for $\Omega_f$ as in Proposition \ref{period}, takes values in $V_{k,k}(E)^*$ for some number field $E$.
\end{itemize}
\end{mprop}
\begin{longversion}
\begin{mrem}For the reader that is unfamiliar with modular symbols, it helps to keep the following mantra in mind. Bianchi forms are inherently analytic in nature, whilst the definition of modular symbols above is purely algebraic. The above proposition allows us to attach a modular symbol to a Bianchi eigenform in a way that retains the data of its Hecke eigenvalues, and it is often easier to study such systems of eigenvalues via this algebraic interpretation, particularly for the purposes of $p$-adic variation.
\end{mrem}
\end{longversion}

\subsection{Modular symbols and $L$-values}
There is an explicit link between critical $L$-values and modular symbols. Recall the notation from Section \ref{notation}. The following is a special case of Theorem 2.11 of \cite{Wil17}. 
\begin{mthm}Let $f \in S_{k,k}(\Omega_0(\n))$ be a cuspidal Bianchi form with associated modular symbol $\phi_f$. Let $\chi$ be a finite order Hecke character of conductor $(c)$ and let $0\leq r\leq k$. Then 
\begin{equation}\label{integralformula}
\Lambda(f,\chi,r+1) = \left[\frac{(-1)^{k}2\Omega_f}{D^{r+1}\tau(\chi^{-1})|c|^r \unitsize}\right]\sum_{a \in (\roi_F/c)^\times} \chi(a)C_{r}(a/c),\end{equation}
where
\[C_{r}(a/c) \defeq \phi_{f}\{a/c - \infty\}((cz+a)^r(\overline{c}\zbar + \overline{a})^r).\]
Here $\tau(\chi^{-1})$ is the Gauss sum attached to $\chi^{-1}$ (see \cite{Wil17}, Section 1.2.3).
\end{mthm}
\begin{longversion}
\begin{mrem}There is a neater version of this theorem in \cite{BW16} in the case of general number fields using the theory of automorphic cycles. We use this version as it is considerably more explicit.
\end{mrem}
\end{longversion}

\subsection{Overconvergent modular symbols}
Part (ii) of Proposition \ref{attach modular symbol} allows us to see the modular symbols $\phi_{f}$ as having values in $V_{k,k}(L)^*$ for a sufficiently large $p$-adic field $L$. For suitable level groups, one can then replace this space of polynomials with a space of $p$-adic distributions to obtain \emph{overconvergent modular symbols}.
\begin{mdef}Let $\mathcal{A}(L)$ denote the space of locally analytic functions on $\roi_F\otimes_{\Z}\Zp$ defined over $L$. We equip this space with a weight $(k,k)$-action of the semigroup
\[\Sigma_0(p) \defeq \left\{\matr \in M_2(\roi_F\otimes_{\Z}\Zp): p|c, a\in(\roi_F\otimes_{\Z}\Zp)^\times, ad-bc \neq 0\right\}\]
by setting
\[\gamma\cdot \zeta (z) = \frac{(a+cz)^k}{(ad-bc)^{k/2}}\zeta\left(\frac{b+dz}{a+cz}\right).\]
\end{mdef}
\begin{mdef}
Let $\mathcal{D}_{k,k}(L) \defeq \mathrm{Hom}_{\mathrm{cts}}(\mathcal{A}(L),L)$ denote the space of locally analytic distributions on $\roi_F\otimes_{\Z}\Zp$ defined over $L$, equipped with a weight $(k,k)$ right action of $\Sigma_0(p)$ given by $\mu|\gamma(\zeta) = \mu(\gamma\cdot\zeta).$ For $\Gamma\subset \Sigma_0(p)$, define the space of \emph{overconvergent modular symbols of weight $(k,k)$ and level $\Gamma$} to be $\symb_{\Gamma}(\mathcal{D}_{k,k}(L)).$
\end{mdef}

There is a natural map $\mathcal{D}_{k,k}(L) \rightarrow V_{k,k}(L)^*$ given by dualising the inclusion of $V_{k,k}(L)$ into $\mathcal{A}(L)$. This induces a \emph{specialisation map}
\[\rho:\symb_{\Gamma_0(\n)}(\mathcal{D}_{k,k}(L)) \longrightarrow \symb_{\Gamma_0(\n)}(V_{k,k}(L)^*),\]
noting that the source is well-defined since $\Gamma_0(\n)\subset \Sigma_0(p)$.
\begin{mthm}[Control theorem]\label{control theorem}
For each prime $\pri$ above $p$, let $\lambda_{\pri}\in L^\times$. If $v(\lambda_{\pri}) < (k+1)/e_{\pri}$ for all $\pri|p$, then the restriction of the specialisation map
\[\rho:\symb_{\Gamma_0(\n)}(\mathcal{D}_{k,k}(L))^{\{U_{\pri}=\lambda_{\pri}:\pri|p\}} \isorightarrow \symb_{\Gamma_0(\n)}(V_{k,k}(L)^*)^{\{U_{\pri}=\lambda_{\pri}:\pri|p\}}\]
to the simultaneous $\lambda_{\pri}$-eigenspaces of the $U_{\pri}$ operators is an isomorphism. Here recall that $e_{\pri}$ is the ramification index of $\pri|p$.
\end{mthm}
\begin{mdef}If $f \in S_{k,k}(\Omega_0(\n))$ is an eigenform with eigenvalues $\lambda_I$, we say $f$ has \emph{small slope} if $v(\lambda_{\pri})<(k+1)/e_{\pri}$ for all $\pri|p$.
\end{mdef}
Thus if $f$ is small slope, using the above control theorem, we get an associated overconvergent modular symbol $\Psi_f \in \symb_{\Gamma_0(\n)}(\mathcal{D}_{k,k}(L))$ by lifting the corresponding classical modular symbol.

\subsection{The $p$-adic $L$-function of a Bianchi modular form}\label{padic lfunction section}
In \cite{Wil17}, the $p$-adic $L$-function of a Bianchi modular form $f$ is defined as a locally analytic distribution on $\cl_F(p^\infty)$ that interpolates the classical $L$-values of $f$. For our purposes, it is more convenient to use a different description in terms of analytic functions on $\roi_F\otimes_{\Z}\Zp$.
First, we need a slight extension of the results of \cite{Wil17}. Let $\mathfrak{g}$ be any ideal coprime to $(p)$. An explicit study of $\cl_F(\mathfrak{g}p^\infty)$ shows that 
\[\cl_F(\mathfrak{g}p^\infty) = \left[(\roi_F/\mathfrak{g})^\times\times(\roi_F\otimes_{\Z}\Zp)^\times\right]/\roi_F^\times.\]
Let $f \in S_{k,k}(\Omega_0(\n))$ be a small slope Bianchi eigenform with associated overconvergent modular symbol $\Psi_f$. Define a distribution $\mu'_{a\newmod{\mathfrak{g}}}$ on $\{[a]\}\times (\roi_F\otimes_{\Z}\Zp) \subset (\roi_F/\mathfrak{g})^\times\times(\roi_F\otimes_{\Z}\Zp)$, which we see as a copy of $\roi_F\otimes_{\Z}\Zp$, by setting 
\[\mu'_{a\newmod{\mathfrak{g}}} \defeq (g\overline{g})^{k/2}\left[\Psi_f\left|\matrd{1}{b}{0}{g}\right.\right]\{0-\infty\},\]
where $b\in \roi_F$ is some lift of $a\newmod{\mathfrak{g}}$ and $\mathfrak{g} = g\roi_F$. This is easily seen to be independent of the choices of $b$ and $g$. Combining these for all $a \in (\roi_F/\mathfrak{g})^\times$, we get a distribution $\mu_p$ on $(\roi_F/\mathfrak{g})^\times \times(\roi_F\otimes_{\Z}\Zp)$. Restricting to units in the second variable, and then restricting to functions that are invariant under $\roi_F^\times$, we obtain a distribution on $\cl_F(\mathfrak{g}p^\infty)$ which in an abuse of notation we will also call $\mu_p$.

\begin{mdef}\label{padiclfunctiondefinition}
For the rest of the paper, fix a choice of $p$-adic logarithm
\[\log_p : \Cp^\times \longrightarrow \Cp\]
with $\log_p(p) = 0$. Also, for $\pri|p$, let $r_{\pri}$ denote the smallest positive integer such that the usual $p$-adic exponential map converges on $\pri^{r_{\pri}}\roi_{\pri}$. (Note that $r_{\pri} = 1$ unless $p = 2$ or $p$ is ramified).
\begin{itemize}
\item[(i)]
Let $\pri|p$ and $s \in \roi_{\pri},$ let $\langle\cdot\rangle : \roi_{\pri}^\times \rightarrow 1+{\pri}^{r_{\pri}}\roi_{\pri}$ denote projection, and define
\begin{align*}
\langle\cdot\rangle^s : \roi_{\pri}^\times &\longrightarrow \Cp,\\
z &\longmapsto \mathrm{exp}(s\cdot\log_p(\langle z \rangle)),
\end{align*}
which is well-defined by the definition of $r_{\pri}$.
\item[(ii)] Similarly, for $\mathbf{s} = (s_{\pri})_{\pri|p} \in\roi_F\otimes_{\Z}\Zp  \cong \prod_{\pri|p}\roi_{\pri}$ define 
\[\langle\cdot\rangle^{\mathbf{s}} \defeq \prod_{\pri|p}\langle\cdot\rangle^{s_{\pri}} : (\roi_F\otimes_{\Z}\Zp)^\times \longrightarrow \Cp,\]
and note that this is invariant under $\roi_F^\times$, so that it induces a map $\langle\cdot\rangle^{\mathbf{s}} : \cl_F(p^\infty) \rightarrow \Cp$. 
\item[(iii)] 
Let $\chi$ be a finite order Hecke character of conductor $\mathfrak{g}\ff$, where $\mathfrak{g}$ is coprime to $(p)$ and $\ff|p^\infty$ (so that $\chi$ can naturally be seen as a finite order character on $\cl_F(\mathfrak{g}p^\infty)$). Define an analytic function on $\roi_F\otimes_{\Z}\Zp$ by 
\[L_p(f,\chi,\mathbf{s}) = \int_{\cl_F(\mathfrak{g}p^\infty)}\langle\mathbf{z}_p\rangle^{\mathbf{s}}\chi(\mathbf{z})d\mu_p(\mathbf{z}),\]
where $\mathbf{s} = (s_{\pri})_{\pri|p} \in \roi_F\otimes_{\Z}\Zp$, and $\mathbf{z}_p =(z_{\pri})_{\pri|p}$ is the projection of $\mathbf{z}$ to $\cl_F(p^\infty)$.
\item[(iv)] Let $\tmp: \roi_{\pri}^\times \rightarrow (\roi_{\pri}/\pri^{r_{\pri}})^\times \subset \roi_{\pri}^\times$ denote the Teichm\"{u}ller character at $\pri$, so that for $x\in\roi_{\pri}^\times,$ we have $x = \tmp(x)\langle x \rangle.$ Also let $\tm \defeq \prod_{\pri|p}\tmp$ be the corresponding character of $(\roi_F\otimes_{\Z}\Zp)^\times.$
\end{itemize}
\end{mdef}
Translated into this setting, the main result of \cite{Wil17} was the following:
\begin{mthm}\label{existence of p-adic L-function}
For a finite order Hecke character $\chi$ with conductor $\mathfrak{g}\ff = (c)$, where $\ff|p^\infty$, and for $0\leq r\leq k$, we have
\[L_p(f,\chi\tm^r,\mathbf{r}) = \left(\prod_{\pri|p}Z_{\pri}(\chi,r)\right)\left[\frac{D^{r+1}\tau(\chi^{-1})|c|^{r} \unitsize}{(-1)^{k}2\lambda_{\ff}\Omega_{f}}\right]\Lambda(f,\chi,r+1),\]
where $\mathbf{r} = (r)_{\pri|p}$ and
\[Z_{\pri}(\chi,r) \defeq 1-\frac{\chi(\pri)N(\pri)^{r}}{\lambda_{\pri}}
\]
recalling that $\chi(\pri)$ is defined to be $\chi(\pi_{\pri})$ (for $\pi_{\pri} \in \roi_{\pri}$ a uniformiser) if $\chi$ is unramified at $\pri$ and $0$ otherwise (so that $Z_{\pri} = 1$ if $\chi$ ramifies at $\pri$). We call $L_p$ the \emph{$p$-adic $L$-function of $f$}.
\end{mthm}
\begin{proof}\emph{(Sketch).} For $\mathfrak{g} = 1$, this is proved entirely in \cite{Wil17}, so we indicate briefly how the more general result follows. Letting $b$ be a unit$\newmod{\mathfrak{g}\ff}$ and using the notation and arguments \emph{op.\ cit.,} Section 7.1, we see that 
\begin{align*}\mu_p(P^{r,r}_{b,\mathfrak{g}\ff}) &= (g\overline{g})^{k/2}\Psi_f\left|\matrd{1}{b}{0}{g}\{0-\infty\}(P^{r,r}_{b,\ff})\right.\\
&= \lambda_{\ff}^{-1}(g\alpha)^r(\overline{g\alpha})^r c_{r,r}\left(\frac{b}{g\alpha}\right),
\end{align*}
where $\ff = (\alpha)$. The rest of the proof proceeds as in \cite{Wil17}.
\end{proof}
\begin{mrems}
\begin{itemize}
\item[(i)] In \cite{Wil17}, this corresponds to the evaluation of $\mu_p$ at the character $\varphi_{p-\mathrm{fin}}$, where $\varphi = \chi|\cdot|^r$. We need the Teichm\"{u}ller character since $\varphi_{p-\mathrm{fin}}(\mathbf{z}) = \chi(\mathbf{z})[\langle \mathbf{z}_p\rangle\tm(\mathbf{z}_p)]^r.$
\item[(ii)]
There is a slight error in \cite{Wil17}, where the term $Z_{\pri}$ is incorrect in the case where $\pri\nmid \ff$. This stems from a subtlety in the definition of Gauss sums. \begin{shortversion}In particular, the proof quotes results from \cite{BW16}, and the Gauss sums in that paper are normalised differently. The difference is the factor of $\varphi_{p-\mathrm{fin}}$.\end{shortversion}
\begin{longversion}In particular, the proof quotes the corresponding (more general) result in \cite{BW16}, which uses an adelic Gauss sum defined in terms of Deligne's epsilon factors, whilst in \cite{Wil17}, the Gauss sum is defined in global terms, and it turns out that these are equal only after multiplying by a suitable scaling factor. To illustrate this, consider the following example over $\Q$; let $\chi$ be a Dirichlet character of conductor $pq$. Then the `global' Gauss sum is defined as $\tau_g(\chi) \defeq \sum_{a \newmod{pq}}\chi(a)e^{2\pi i a/pq}$, whilst the `adelic' one is defined as $\tau_a(\chi) \defeq \sum_{a \newmod{pq}}\chi(a)e^{2\pi i (a/p + a/q)} = \chi(p+q)^{-1}\tau_g(\chi).$ In the class number 1 case over imaginary quadratic fields, the result on the exceptional factor is to precisely cancel the factor of $\varphi_{p-\mathrm{fin}}(\pi_{\pri})$ from $Z_{\pri}$, leaving the formulation above. (Regrettably, there is another discrepancy between \cite{Wil17} and \cite{BW16} in that the Gauss sum of a character in the former relates to the Gauss sum of the inverse character in the latter. Throughout this paper we use the normalisations of the former).
\item[(iii)]
For $\mathfrak{g} = 1$, the interpolation property proved in \cite{Wil17} is more general, allowing characters of infinity type $(q,r)$ where both $0 \leq q,r \leq k$. For our purposes it suffices to only consider characters of the form $\chi|\cdot|^r$ as in the theorem above.
\end{longversion}
\end{itemize}
\end{mrems}

%
%
\section{Modular symbols on the tree}
\label{modular symbols on the tree}

\subsection{Harmonic cocycles attached to Bianchi modular forms}
Working with class number 1 allows us to work classically with forms on the upper half-space, for which we'll need a classical version of the (adelic) group $\Omega$ introduced in Definition \ref{group omega}.
\begin{mdef}\label{group gamma}
Define 
\[\Gamma \defeq \mathrm{PGL}_2(F)\cap \Omega.\]
The group $\Gamma$ has an explicit description as the subgroup of $\PGLt(\roi_F[\pri^{-1}])$ with lower left entry contained in $\pri^{r}\m$, for some (possibly negative) $r\in\Z$, whose determinant is exactly divisible by an even power of $\pri$.
\end{mdef}
Let $f \in S_{k,k}(\Omega_0(\n))$ be a Bianchi modular form that is new at $\pri$, and let $\f$ be the associated form on the tree $\treep$. To $\f$ we can associate a `modular symbol' $\kappa_{\f}$ on $\treep$ by setting
\[\kappa_{\f}\{r-s\}(e,P) \defeq \phi_{\f_{e}}\{r-s\}(P),\]
for $e \in \edges(\treep), r,s \in \Proj(F_{\pri})$ and $P\in V_{k,k}(L)$, for $L$ a suitably large finite extension of $\Qp$.
\begin{mprop}\label{modular symbol tree properties}
The symbol $\kappa_{\f}$ is harmonic on $\edges(\treep)$, linear in $r, s$ and $P$ and $\Gamma$-invariant in the sense that
\[\kappa_{\f}\{\gamma r - \gamma s\}(\gamma e, \gamma\cdot P) = \kappa_{\f}\{r-s\}(e,P)\]
for all $\gamma \in \Gamma$.
\end{mprop}
\begin{proof}(Sketch). These properties all follow from the corresponding properties of $\f$ and $\phi_{\f_{e}}$. The only one that is not obvious is the $\Gamma$-invariance property, which follows from a lengthy technical calculation, of which we give a sketch. Firstly, one can obtain a more classical definition of Bianchi modular forms on $\treep$ as functions $\f : \edges(\treep)\times\uhs \longrightarrow V_{2k+2}(\C)$ such that
\[\f(\gamma e, \gamma(z,t))\left[\binomc{X}{Y}\right] = \f(e, (z,t))\left[j(\gamma,(z,t))\binomc{X}{Y}\right],\]
for all $\gamma = \smallmatrd{a}{b}{c}{d} \in \Gamma$, where we consider $V_{2k+2}(\C)$ instead as homogeneous polynomials in two variables $X$ and $Y$ and
\[j(\gamma,(z,t)) \defeq \matrd{c}{0}{0}{\overline{c}}\matrd{z}{-t}{t}{\overline{z}} + \matrd{d}{0}{0}{\overline{d}}.\]
This is precisely the transformation property that classical Bianchi modular forms satisfy under the corresponding level group. Then one examines the definition of the modular symbol $\phi_{\f_e}$, via a $V_{k,k}(\C)^*$-valued harmonic differential $\delta_{\f_e}$ on $\uhs$ attached to $\f_e$. In \cite{Hid94}, Chapter 2.5, this differential is shown to be invariant under $\Gamma_e$. Define a differential on $\edges\times\uhs$ by $\delta(e,(z,t)) \defeq \delta_{\f_e}(z,t)$. Using exactly the same methods as in \cite{Hid94}, and the classical transformation property above, this is seen to be invariant under $\Gamma$, in the sense that $\delta(\gamma e, \gamma(z,t)) = \delta(e,(z,t))$ for all $\gamma \in \Gamma$. Since by definition $\kappa_{\f}\{r-s\}(e,P)$ is the integral of $\delta(e,(z,t))$ from $r$ to $s$ in $\uhs$, we then obtain the result by a simple transformation.
\end{proof}

\begin{mdef}
We call such a function a \emph{$\Gamma$-invariant harmonic cocycle on $\edges(\treep)$ with values in $\mathrm{Hom}(\Delta_0,V_{k,k}(L)^*).$} We denote the space of such functions by \linebreak[4]$C^{\mathrm{har}}(\mathrm{Hom}(\Delta_0,V_{k,k}(L)^*)^{\Gamma}.$
\end{mdef}
We have a map
\[\rho_{e_*} : C^{\mathrm{har}}(\mathrm{Hom}(\Delta_0,V_{k,k}(L)^*)^{\Gamma} \longrightarrow \symb_{\Gamma_0(\n)}(V_{k,k}(L)^*)\]
given by restriction to the standard edge $e_*$. We can implicitly define an inverse to this map by using the tree. In fact, we can write this down explicitly without needing to refer to the tree at all. Indeed, note that any \emph{even} edge $e$ can be written as $e = \gamma e_*$ for some $\gamma \in \Gamma$. Let $\phi_{f} \in \symb_{\Gamma_0(\n)}(V_{k,k}(L)^*)$ be a modular symbol. If $\kappa_{\f}$ is a harmonic cocycle such that $\rho_{e_*}(\kappa_{\f}) = \phi_{f},$ then by $\Gamma$-invariance we have
\[\kappa_{\f}\{r - s\}(e, P) = \phi_{f}|\gamma^{-1}\{r-s\}(P).\]
It remains to determine the value of $\kappa_{\f}$ at odd edges, which we will do in the next subsection using the \emph{Atkin--Lehner} operators to describe the action of a larger group $\widetilde{\Gamma}.$
\begin{mrem}
We have brushed over a slight subtlety here; we should really have defined $\kappa_{\f}$ in terms of the complex symbols $\widetilde{\phi}_{\f_e}$, rather than the algebraic analogues $\phi_{\f_e}$, each of which depends on a choice of period $\Omega_{\f_e}$. However, the invariance property above shows that we can use the same choice $\Omega_{\f_e} = \Omega_f$ at each edge to obtain something algebraic on the tree.
\end{mrem}

\subsection{Atkin--Lehner operators and the action of $\widetilde{\Gamma}$}

Let $\widetilde{\Gamma} \defeq \PGLt(F)\cap\Omegat$. Analogously to $\Gamma$, we see $\Gammat$ as the subgroup of $\PGLt(\roi_F[\pri^{-1}])$ with lower left entry contained in $\pri^{r}\m$, for some (possibly negative) $r\in\Z$.

\begin{mdef}Let $\alpha$ be any element of $\widetilde{\Gamma}\backslash \Gamma$ normalising $\Gamma_0(\n)$. The \emph{Atkin--Lehner operator} $W_{\pri}$ on $\symb_{\Gamma_0(\n)}(V_{k,k}(L)^*)$ is defined by 
\[\phi|W_{\pri} \defeq \phi|\alpha.\]
If $f$ (and hence $\phi_f$) is an eigenform for the Hecke operators, then it is also an eigenform for $W_{\pri}$. Denote the eigenvalue by $-\omega$, which will be $1$ or $-1$.
\end{mdef}
\begin{mprop}\label{transform}
For $\gamma\in\GLt(F)$, define $|\gamma| \defeq \mathrm{ord}_{\pri}(\det(\gamma))$, so that in particular $|\gamma| \newmod{2}$ is well-defined on $\PGLt(F)$. Then we have
\[\kappa_{\f}\{\gamma r - \gamma s\}(\gamma e, \gamma\cdot P) = \omega^{|\gamma|}\kappa_{\f}\{r-s\}(e,P)\]
for all $\gamma \in \widetilde{\Gamma}.$
\end{mprop}
\begin{proof}
If $\gamma \in \Gamma$, this is the invariance statement above. If $\gamma \notin \Gamma$, then $\gamma = \gamma'\alpha$ for some $\gamma' \in \Gamma$, and then the $\Gamma$-invariance means we're reduced to proving that 
\[\kappa_f\{\alpha r -\alpha s\}(\alpha e,\alpha \cdot P) = \omega \kappa_f\{r-s\}(e,P).\]
\begin{longversion}To see this, suppose first that $e$ is even, so that $e = \gamma e_*$ for some $\gamma \in \Gamma$. Then
\begin{align*}\omega\kappa_f\{r-s\}(e,P) &= \omega\kappa_f\{\gamma^{-1}r-\gamma^{-1}s\}(e_*,\gamma^{-1}\cdot P)\\
&= \omega\phi_f\{\gamma^{-1}r-\gamma^{-1}s\}(\gamma^{-1}\cdot P)\\
&= -(\phi_f|W_{\pri})\{\gamma^{-1}r-\gamma^{-1}s\}(\gamma^{-1}\cdot P) = -\phi_f\{\alpha\gamma^{-1}r-\alpha\gamma^{-1}s\}(\alpha\gamma^{-1}\cdot P)\\
&= -\kappa_f\{\alpha\gamma^{-1}r-\alpha\gamma^{-1}s\}(e_*,\alpha\gamma^{-1}\cdot P)\\
&= \kappa_f\{\alpha\gamma^{-1}r-\alpha\gamma^{-1}s\}(\overline{e_*},\alpha\gamma^{-1}\cdot P)\\
&= \kappa_f\{\gamma\alpha\gamma^{-1}r-\gamma\alpha\gamma^{-1}s\}(\overline{e},\gamma\alpha\gamma^{-1}\cdot P)\\
&= \kappa_f\{\gamma\alpha\gamma^{-1}r-\gamma\alpha\gamma^{-1}s\}(\gamma\alpha\gamma^{-1}e,\gamma\alpha\gamma^{-1}\cdot P).
\end{align*}
But $\gamma\alpha\gamma^{-1} \in \Gammat\backslash\Gamma$, so it is of the form $\gamma'\alpha$ for some $\gamma' \in \Gamma$.  The result follows since we can cancel $\gamma'$ by $\Gamma$ invariance \end{longversion}
\begin{shortversion}This is a simple explicit calculation \end{shortversion}(compare \cite{Dar01}, Lemma 1.4).
\begin{longversion}
$\lb$
The case where $e$ is odd is similar, using the existence of some $\gamma$ such that $e = \gamma\alpha e_*$. This concludes the proof.
\end{longversion}
\end{proof}

\subsection{Hecke operators}
There are natural actions of Hecke operators away from $\n$ on the forms and symbols on the tree, and all the associations we've made are equivariant with respect to these operators. In particular, recall Definition \ref{delta} where for each ideal $I$, we picked a finite set of representatives $\delta_j^I$ for the double coset corresponding to $T_I$.
\begin{mdef}
Define the Hecke operator $T_I$ on $C^{\mathrm{har}}(\mathrm{Hom}(\Delta_0,V_{k,k}(L)^*)^{\Gamma}$ by setting $\kappa|T_I \defeq N(I)^{k/2}\sum_j\kappa|\delta_j^I,$ or more concretely,
\[(\kappa|T_I)\{r-s\}(e,P) \defeq N(I)^{k/2}\sum_j\kappa\{\delta_j^I r - \delta_j^Is\}(\delta_j^Ie, \delta_j^I\cdot P).\]
\end{mdef}
\begin{mprop}
The association of $\kappa_{\f}$ to $f$ is equivariant with respect to the Hecke operators $T_I$ for $I$ coprime to $\n$.
\end{mprop}
\begin{proof}
This follows immediately from the definitions and the fact that, if $I$ is coprime to $\n$, we have $\delta_j^I e = e$ for all $j$ and for any edge $e$.
\end{proof}

\subsection{Hecke eigenvalues of newforms}
As a corollary of the above results, we recover the following classical result, which also follows from the theory of automorphic forms and the structure of Iwahori--Hecke algebras. The authors thank Aurel Page for pointing out the previous existence of such a result, which we were not able to find in the literature.

\begin{mcor}\label{eigenvalues of newforms}
Let $f\in S_{k,k}(\Omega_0(\n))$ be an eigenform, with $\pri|\n$. If $f$ is new at $\pri$, then the Hecke eigenvalue at $\pri$ is $\omega N(\pri)^{k/2}$, where $-\omega$ is the eigenvalue at $f$ of the Atkin--Lehner operator $W_{\pri}$.
\end{mcor}
\begin{proof}
We've already shown that $f$ is new at $\pri$ if and only if the associated form $\f$ on the tree is harmonic; suppose this is the case. It suffices to prove the result on modular symbols. Recall that $\pi$ is a generator of $\pri$, and for $j \in \roi_F$, let $\delta_j \defeq \smallmatrd{1}{j}{0}{\pi}$. Then by definition, for any $r,s \in \Proj(F)$ and $P\in V_{k,k}(L),$
\begin{align*}(\phi_f|U_{\pri})\{r-s\}(P) &= N(\pri)^{k/2}\sum_{j \newmod{\pri}}\kappa_{\f}\{\delta_j r - \delta_j s\}(e_*, \delta_j\cdot P)\\
&= \omega N(\pri)^{k/2}\sum_{j \newmod{\pri}}\kappa_{\f}\{r-s\}(\delta_j^{-1}e_*, P),
\end{align*}
using the transformation property above. But by inspection, the set $\{\delta_j^{-1}e_*: j \newmod{\pri}\}$ is precisely the set $\{e \in \edges(\treep): t(e) = v_*, e\neq \overline{e_*}\}.$ 
\begin{longversion}Hence we get
\begin{align*}(\phi_f|U_{\pri})\{r-s\}(P) &= \omega N(\pri)^{k/2}\sum_{\substack{t(e) = v_*\\e \neq \overline{e_*}}} \kappa_f\{r-s\}(e,P)\\
&= -\omega N(\pri)^{k/2}\kappa_f\{r-s\}(\overline{e_*},P)\\
&= \omega N(\pri)^{k/2}\kappa_f\{r-s\}(e_*, P),
\end{align*}
where the last two equalities follow from harmonicity. This completes the proof.
\end{longversion}
\begin{shortversion}Hence an easy calculation using harmonicity shows that
\[(\phi_f|U_{\pri})\{r-s\}(P) = \omega N(\pri)^{k/2}\kappa_f\{r-s\}(e_*,P),\]
which completes the proof.
\end{shortversion}
\end{proof}
Whilst we have presented this in only the generality in which we require it, an argument of this type generalises very easily to the case of arbitrary number fields and arbitrary class number.
\begin{mcor}\label{exceptional zero existence}
Let $f$ as above be new at $\pri = (\pi)$. Suppose that $f$ has small slope at all primes above $p$, and let $\chi$ be a finite order Hecke character for $F$. Then the critical value $L_p(f,\chi\tm^{r},\mathbf{r})$ of the $p$-adic $L$-function has an exceptional zero at $\pri$ precisely when $k$ is even, $\chi(\pi) = \omega$ and $r = k/2$.
\end{mcor}
\begin{proof}
The exceptional factor $Z_{\pri}(\chi,k/2)$ of Theorem \ref{existence of p-adic L-function} is zero if and only if these conditions are satisfied.
\end{proof}

%
%
\section{Distributions attached to Bianchi modular forms}
In the rational set-up, overconvergent modular symbols naturally take values in the space of locally analytic distributions on $\Zp$, and Darmon and Orton use the tree to extend these to distributions on $\Proj(\Qp)$. In the Bianchi case, the values are distributions on $\roi_F\otimes_{\Z}\Zp$. In this section we show how to extend such distributions to be `projective at $\pri$'. To do so in general, we'll need to work with the other primes above $p$, but to ease the exposition, we first treat the case where there is a single prime above $p$.

\subsection{Distributions on $\Proj_{\pri}$: the case $p$ inert or ramified}
Suppose $p$ is inert or ramified in $F$. To ensure consistent notation with the split case, let
\[\Proj_{\pri} \defeq \ProjF.\]

\subsubsection{Polynomial distributions on $\Proj_{\pri}$}
\begin{longversion}
Recall Proposition \ref{basis of opens}, which said that edges of the tree give a basis of open sets in $\Proj_{\pri}$ by associating to $e \in \edges(\treep)$ the open set $U(e)$. Let $f\in S_{k,k}(\Omega_0(\n))^{\pri-\mathrm{new}}$ be an eigenform with associated form $\f$ and modular symbol $\kappa_{\f}$ on $\treep$. To this, we associate a system of locally polynomial distributions.
\begin{mdef}
Let $\mathcal{V}_{k,k}(\Proj_{\pri},L)$ be the space of functions $\Proj_{\pri} \rightarrow L$ that are locally polynomial of degree at most $k$ in each variable.
\end{mdef}
\begin{mdef}
For each pair $r,s \in \Proj(F)$, define a locally polynomial distribution $\mu_{\f}^{\mathrm{poly}}\{r-s\}$ on $\Proj_{\pri}$ by setting
\[\int_{U(e)}P(t)d\mu_{\f}^{\mathrm{poly}}\{r-s\}(t) \defeq \kappa_{\f}\{r-s\}(e,P),\]
where $P$ is polynomial of degree at most $k$ in each variable on $U(e)$, and extending linearly.
\end{mdef}

Recall that $-\omega$ is the eigenvalue of $f$ at the Atkin--Lehner operator $W_{\pri}$. The transformation property of Proposition \ref{transform} above immediately gives:
\begin{mlem}\label{poly distribution}
Let $e = \gamma e_*$ with $\gamma \in \Gammat$. Then 
\[\int_{U(e)}P(t) d\mu_{\f}^{\mathrm{poly}}\{r-s\}(t) = \omega^{|\gamma|}\phi_f\{\gamma^{-1}r -\gamma^{-1}s\}(\gamma^{-1}\cdot P).\]
\end{mlem}
In particular, we can express the distribution solely in terms of $\phi_f$.

\begin{mrems}
\begin{itemize}
\item[(i)]This distribution is well-defined by the harmonicity property of $\kappa_{\f}$. In particular, constructions of the sort above are why we need the harmonicity condition.
\item[(ii)] This construction, which is directly analogous to the one in \cite{Ort04} for classical modular forms, is included mainly to motivate future constructions of locally analytic distributions. We're really seeing the modular symbol attached to $\f_e$ to take values in `polynomials on $U(e)$.'
\end{itemize}
\end{mrems}

\end{longversion}
\begin{shortversion}By Proposition \ref{basis of opens}, a basis of open sets in $\Proj_{\pri}$ is given by $\{U(e): e \in \edges(\treep)\}.$ Let $f\in S_{k,k}(\Omega_0(\n))^{\pri-\mathrm{new}}$ be an eigenform with associated form $\f$ and modular symbol $\kappa_{\f}$ on $\treep$. Let $\mathcal{V}_{k,k}(\Proj_{\pri})$ be the space of functions $\Proj_{\pri} \rightarrow L$ that are locally polynomial of degree at most $k$ in each variable.
\begin{mdef}
For each pair $r,s \in \Proj(F)$, define a locally polynomial distribution $\mu_{\f}^{\mathrm{poly}}\{r-s\}$ on $\Proj_{\pri}$ by setting
\[\int_{U(e)}P(t)d\mu_{\f}^{\mathrm{poly}}\{r-s\}(t) \defeq \kappa_{\f}\{r-s\}(e,P),\]
where $P$ is polynomial on $U(e)$ of degree at most $k$ in each variable, and extending linearly. (Note that this is well-defined by harmonicity).
\end{mdef}
Recall that $-\omega$ is the eigenvalue of $f$ at the Atkin--Lehner operator $W_{\pri}$. The transformation property of Proposition \ref{transform} above then shows that if $e = \gamma e_*$ with $\gamma \in \Gammat$, then 
\begin{equation}\label{poly distribution}\int_{U(e)}P(t) d\mu_{\f}^{\mathrm{poly}}\{r-s\}(t) = \omega^{|\gamma|}\phi_f\{\gamma^{-1}r -\gamma^{-1}s\}(\gamma^{-1}\cdot P).
\end{equation}
In particular, we can express the distribution solely in terms of $\phi_f$.

\end{shortversion}

\subsubsection{Locally analytic distributions on $\Proj_{\pri}$}
By Corollary \ref{eigenvalues of newforms}, $f$ has eigenvalue $\omega N(\pri)^{k/2} = \omega p^{k/e_{\pri}}$ at $\pri$, so $f$ has small slope at $p$. Hence we can use Theorem \ref{control theorem} to give a unique overconvergent lift $\Psi_f\in\symb_{\Gamma}(\mathcal{D}_{k,k}(L))$ of $\phi_f$.
\begin{mdef}
Let $\mathcal{A}_{k,k}(\Proj_{\pri},L)$ denote the space of functions $\Proj_{\pri} = \ProjF \rightarrow L$ that are locally analytic except perhaps for a pole at $\infty$ of order at most $k$.
\end{mdef}
Recall the action of $\PGLt(F_{\pri})$ on $\Proj_{\pri}$ by $\smallmatrd{a}{b}{c}{d}\cdot x = (b+dx)/(a+cx).$ We get an action of $\PGLt(F_{\pri})$ on $\mathcal{A}_{k,k}(\Proj_{\pri},L)$ by
\[\gamma\cdot \zeta(x) \defeq \frac{(a+cx)^k(\overline{a}+\overline{cz})^k}{(ad-bc)^{k/2}(\overline{ad}-\overline{bc})^{k/2}}\zeta(\gamma\cdot x),\]
directly extending the action of $\Sigma_0(p)$ on locally analytic functions on $\roi_{\pri} = \roi_F\otimes_{\Z}\Zp$ from earlier.\\
\\
Note that if $e = \gamma e_*$, then $\gamma^{-1}$ carries the set $\roi_{\pri}$ to $U(e) = \gamma^{-1}(\roi_{\pri}) \subset \Proj_{\pri}$. Hence if $\zeta$ is supported on $U(e)$, then $\gamma^{-1}\cdot \zeta$ is supported on $\roi_{\pri}$.
\begin{longversion} This, along with Lemma \ref{poly distribution}, motivates the following.
\end{longversion}
\begin{shortversion}
This, along with equation (\ref{poly distribution}), motivates the following.
\end{shortversion}
\begin{mdef}Let $e = \gamma e_* \in \edges(\treep)$.
\begin{itemize}
\item[(i)]Let $\mathcal{A}_{k,k}(U(e),L)$ be the subspace of functions $\zeta \in \mathcal{A}_{k,k}(\Proj_{\pri},L)$ that are supported on $U(e)$, with the natural weight $k$ action of $\gamma\Sigma_0(p)\gamma^{-1}$, and let $\mathcal{D}_{k,k}(U(e),L)$ be its continuous dual.
\item[(ii)] Define the \emph{overconvergent modular symbol at $e$ associated to $\f$}, denoted by
\begin{longversion}
\[\Psi_{\f_e} \in \symb_{\Gamma_e}(\mathcal{D}_{k,k}(U(e),L)),\]
by
\end{longversion}
\begin{shortversion}
$\Psi_{\f_e}$, to be the element of $\symb_{\Gamma_e}(\mathcal{D}_{k,k}(U(e),L))$ defined by
\end{shortversion}
\[\Psi_{\f_e}\{r-s\}(\zeta) \defeq \omega^{|\gamma|}\Psi_f\{\gamma^{-1}r-\gamma^{-1}s\}(\gamma^{-1}\cdot \zeta).\]
Here $\Gamma_e \defeq \mathrm{Stab}_{\Gamma}(e) = \gamma\Gamma_0(\n)\gamma^{-1}$.
\end{itemize}
\end{mdef}
\begin{mprop}
For each pair $r,s \in \Proj(F)$, we may define a distribution $\mu_{\f}\{r-s\}$ on $\mathcal{A}_{k,k}(\Proj_{\pri},L)$ attached to $\f$ by setting
\[\int_{U(e)} \zeta(t)d\mu_{\f}\{r-s\}(t) \defeq \Psi_{\f_e}\{r-s\}(\zeta)\]
and extending linearly.
\end{mprop}

\begin{longversion}
\begin{proof}
Suppose $p$ is inert (the case $p$ ramified is essentially identical). We need only show that $\mu_{\f}\{r-s\}$ is distributive over open sets. Let $e = \gamma e_*\in \edges(\treep)$ be an edge, with source $v \in \vertices(\tree)$, and let $\overline{e}, e_1,...,e_{p^2}$ be the edges with target $v$, so that $U(e) = \coprod _{i=1}^{p^2}U(e_i).$ We need to show that, for $\zeta$ locally analytic on $U(e)$, we have
\[\Psi_{\f_e}\{r-s\}(\zeta) = \sum_{i=1}^{p^2}\Psi_{\f_{e_i}}\{r-s\}(\zeta|_{U(e_i)}).\]
For each $j \newmod{p}$, recall that $\delta_j = \smallmatrd{1}{j}{0}{p}$. By direct inspection of the tree, we see that
\[\{e_i: 1\leq i \leq p^2\} = \{\delta_j^{-1} \gamma e_*: j \newmod{p}\} = \{\gamma \delta_j^{-1} e_*: j\newmod{p}\}.\]
(We need not have $\delta_j^{-1} \gamma e_* = \gamma \delta_j^{-1} e_*$ for individual edges, but we do get an equality of sets). Hence
\begin{align*}\sum_{i=1}^{p^2}\Psi_{\f_{e_i}}\{r-s\}(\zeta|_{U(e_i)}) &= \omega^{1+|\gamma|}\sum_{j \newmod{p}}\Psi_{f}\{\delta_j\gamma^{-1} r- \delta_j\gamma^{-1} s\}(\delta_j\gamma^{-1}\cdot \zeta)\\
&=  \omega^{1+|\gamma|}N(\pri)^{-k/2}(\Psi_{f}|U_{\pri})\{\gamma^{-1}r-\gamma^{-1}s\}(\gamma^{-1}\cdot\zeta)\\
&= \omega^{|\gamma|}\Psi_f\{\gamma^{-1}r-\gamma^{-1}s\}(\gamma^{-1}\cdot\zeta) = \Psi_{\f_e}\{r-s\}(\zeta),
\end{align*}
as required. Here the penultimate equality follows from Corollary \ref{eigenvalues of newforms}.
\end{longversion}

\begin{shortversion}
\begin{proof}
Suppose $p$ is inert (the case $p$ ramified is essentially identical). We need only show that $\mu_{\f}\{r-s\}$ is distributive over open sets. Let $e = \gamma e_*\in \edges(\treep)$ be an edge, with source $v \in \vertices(\tree)$, and let $\overline{e}, e_1,...,e_{p^2}$ be the edges with target $v$, so that $U(e) = \coprod _{i=1}^{p^2}U(e_i).$ We need to show that, for $\zeta$ locally analytic on $U(e)$, we have
\[\Psi_{\f_e}\{r-s\}(\zeta) = \sum_{i=1}^{p^2}\Psi_{\f_{e_i}}\{r-s\}(\zeta|_{U(e_i)}).\]
For each $j \newmod{p}$, recall that $\delta_j = \smallmatrd{1}{j}{0}{p}$. By direct inspection of the tree, we see that $\{e_i: 1\leq i \leq p^2\} = \{\delta_j^{-1} \gamma e_*: j \newmod{p}\} = \{\gamma \delta_j^{-1} e_*: j\newmod{p}\}.$ (We need not have $\delta_j^{-1} \gamma e_* = \gamma \delta_j^{-1} e_*$ for individual edges, but we do get an equality of sets). A calculation using harmonicity and the fact that $\Psi_{f} = \omega p^{-k}\Psi_f|U_{\pri}$ (from Corollary \ref{eigenvalues of newforms}) shows that
\[\sum_{i=1}^{p^2}\Psi_{\f_{e_i}}\{r-s\}(\zeta|_{U(e_i)}) = \omega^{|\gamma|}\Psi_f\{\gamma^{-1}r-\gamma^{-1}s\}(\gamma^{-1}\cdot\zeta) = \Psi_{\f_e}\{r-s\}(\zeta),
\]
as required. 
\end{proof}
\end{shortversion}

\subsection{Distributions on $\Proj_{\pri}$: the case $p$ split}
Now suppose $p$ is split as $\pri\pribar$ in $F$. The only additional complication in this case is that to construct $p$-adic $L$-functions, we must also consider behaviour at $\pribar$. Other than this, the construction is closely analogous to the case $p$ inert, and we will omit the details in this case.
\begin{mass}
Henceforth, we will assume that $f$ has small slope at $\pribar$, that is, that $v_p(a_{\pribar}) < k+1$, where $a_{\pribar}$ is the Hecke eigenvalue at $\pribar$. (Note that $f$ automatically has small slope at $\pri$ since it is new, so in this case we obtain a canonical overconvergent symbol $\Psi_f$, as before).
\end{mass}
\begin{mdef}
\begin{itemize}
\item[(i)]Let
\[\Proj_{\pri} \defeq \Proj(F_{\pri})\times \roi_{\pribar} \cong \Proj(\Qp)\times\Zp \supset \roi_F\otimes_{\Z}\Zp.\]
\item[(ii)]Let $\mathcal{A}_k(\ProjF,L)$ be the space of functions $\ProjF \rightarrow L$ that are locally analytic except perhaps for a pole of order at most $k$ at $\infty$, and let $\mathcal{A}_k(\roi_{\pribar},L)$ be the space of locally analytic functions $\roi_{\pribar} \rightarrow L$. Define
\[\mathcal{A}_{k,k}(\Proj_{\pri},L) \defeq \mathcal{A}_k(\ProjF,L)\widehat{\otimes}_{L}\mathcal{A}_k(\roi_{\pribar},L).\]
\end{itemize}
\end{mdef}
Let $U'(e) \defeq U(e)\times \roi_{\pribar} \subset \Proj_{\pri}$. Note that $\Gammat$ preserves $\roi_{\pribar}$, so if $e = \gamma e_*$ for $\gamma \in \Gammat$, then $U'(e) = \gamma^{-1}(\roi_F\otimes_{\Z}\Zp)$. Hence the following definition makes sense.
\begin{mdef}Let $e = \gamma e_* \in \edges(\treep)$. 
\begin{itemize}
\item[(i)]Let $\mathcal{A}_{k,k}(U'(e),L)$ be the subspace of functions $\zeta \in \mathcal{A}_{k,k}(\Proj_{\pri},L)$ that are supported on $U'(e)$ and let $\mathcal{D}_{k,k}(U'(e),L)$ be its continuous dual.
\item[(ii)]Define the \emph{overconvergent modular symbol at $e$ associated to $\f$}, denoted by $\Psi_{\f_e},$ to be the element in $\symb_{\Gamma_e}(\mathcal{D}_{k,k}(U'(e),L))$ defined by
\[\Psi_{\f_e}\{r-s\}(\zeta) \defeq \omega^{|\gamma|}\Psi_f\{\gamma^{-1}r-\gamma^{-1}s\}(\gamma^{-1}\cdot \zeta).\]
\end{itemize}
\end{mdef}
In exactly the same manner as the inert case, we can then show that we may attach a distribution $\mu_{\f}\{r-s\}$ on $\mathcal{A}_{k,k}(\Proj_{\pri},L)$ to $\f$ by defining it on basic open sets as
\[\int_{U'(e)} \zeta(t)d\mu_{\f}\{r-s\}(t) \defeq \Psi_{\f_e}\{r-s\}(\zeta).\]

\subsection{Properties of the distributions}
Regardless of the splitting behaviour of $p$, we've constructed a family of distributions $\mu_{\f}\{r-s\}$ on $\Proj_{\pri}$ attached to $\f$. We have the following easy properties (see \cite{Ort04}, Lemma 3.2).
\begin{mprop}\label{properties of distribution}
\begin{itemize}
\item[(i)] If $P$ is any polynomial function over $L$, then for any $r,s \in \Proj(F)$ we have
\[\int_{\Proj_{\pri}} P(t)d\mu_{\f}\{r-s\}(t) = 0.\]
\item[(ii)]For $\gamma \in \Gammat$, we have
\[\int_{U(\gamma e)}(\gamma\cdot\zeta)(t)d\mu_{\f}\{\gamma r - \gamma s\}(t) =  \omega^{|\gamma|}\int_{U(e)}\zeta(t)d\mu_{\f}\{r - s\}(t).\]
\item[(iii)]The distributions `behave well under Hecke operators' in the sense that for all $I$ coprime to $\n$,
\[N(I)^{k/2}\sum_{j}\int_{U(\delta_{j}^{I} e)}(\delta_{j}^{I}\cdot \zeta)(t)d\mu\{\delta_{j}^{I} r - \delta_{j}^{I} s\}(t) = \lambda_{I}\int_{U(e)}\zeta(t)d\mu_{\f}\{r-s\}(t),\]
keeping the notation of Definition \ref{delta}, and where $\lambda_I$ is the Hecke eigenvalue of $f$ at $I$.

\end{itemize}
\end{mprop}

%
%

\section{Double integrals and cohomology classes}
In this section, we use the machinery constructed above to define double integrals in the manner of Darmon, and then attach two cohomology classes $\lc, \oc \in \HGamma$ to $f$. This closely follows the methods of \cite{Ort04}, which was motivated by the methods of \cite{Dar01}. Darmon defines a \emph{multiplicative} double integral and a class $\mathrm{c}_f$ attached to a weight two modular form $f$, and then constructs classes $\log_p(\mathrm{c}_f)$ and $\mathrm{ord}_p(\mathrm{c}_f)$. In our setting, this multiplicative integral doesn't exist, and we instead define the analogues of $\log_p(\mathrm{c}_f)$ and $\mathrm{ord}_p(\mathrm{c}_f)$ by exploiting their properties as given in \cite{Dar01}.

\subsection{Double integrals}
Let $f \in S_{k,k}(\Omega_0(\n))^{\pri-\mathrm{new}}$ with associated system of distributions $\mu_{\f}\{r-s\}$ on $\Proj_{\pri}$, and recall that we fixed a choice of $p$-adic logarithm $\log_p : \Cp^\times \longrightarrow \Cp$.
\begin{mdef}
 Let $x,y$ in $\uhp_{\pri}$, let $P \in V_{k,k}(\Cp)$ and let $r,s \in \Proj(F)$. Then define
\[\int_x^y\int_r^s (P) \defeq \int_{\Proj_{\pri}}\log_p\left(\frac{t_{\pri}-x}{t_{\pri}-y}\right)P(t)d\mu_{\f}\{r-s\}(t),\]
where $t_{\pri}$ is the projection of $t \in \Proj_{\pri}$ to $\Proj(F_{\pri})$. (Note that this is well-defined since $\uhp_{\pri} \defeq \Proj(\Cp)\backslash\Proj(F_{\pri})$, so that $t_{\pri}-x$ and $t_{\pri}-y$ cannot vanish and their quotient gives an element of $\Cp^\times$).
\end{mdef}
The following (compare \cite{Ort04}, Lemma 4.1) follow from properties of $\log_p$ and from Proposition \ref{properties of distribution}.
\begin{mprop}[Properties of the double integral]\label{properties of double integral}
\begin{itemize}
\item[(i)] The double integral is additive in $x,y$ and $r,s$ and is linear in $P$.
\item[(ii)] For $\gamma\in\Gammat$,
\[\int_{\gamma x}^{\gamma y}\int_{\gamma r}^{\gamma s}(\gamma\cdot P) = \omega^{|\gamma|}\int_x^y\int_r^s (P).\]
\item[(iii)]The double integral `behaves well under Hecke operators' in the sense that for all $I$ coprime to $\n$,
\[N(I)^{k/2}\sum_{j}\int_{\delta_j^I x}^{\delta_k^I y}\int_{\delta_j^I r}^{\delta_j^I s} (\delta_k^I\cdot P) = \lambda_I\int_x^y\int_r^s (P),\]
for notation as in Proposition \ref{properties of distribution} (iii).
\end{itemize}
\end{mprop}

\begin{shortversion}
For ease of notation, let $\Delta_{k,k} \defeq \mathrm{Hom}(\Delta_0,V_{k,k}(\Cp)^*).$ We define cohomology classes $\lc$ and $\oc$ as follows.
\begin{mdef}[Definition of $\lc$] \label{define lc}Let $\tau \in \uhp_{\pri}$, and define a function $\Gamma \rightarrow \Delta_{k,k}$ by
\[\lct(\gamma)\{r-s\}(P) \defeq \int_{\tau}^{\gamma \tau}\int_r^s (P).\]
From the properties of the double integral above, we see that this is a cocycle in $\mathrm{Z}^1(\Gamma,\Delta_{k,k})$ and its class $\lc \in \HGammas$ is independent of $\tau$.
\end{mdef}
\begin{mdef}[Definition of $\oc$] Let $v \in \vertices(\treep)$, and define a function $\Gamma \rightarrow \Delta_{k,k}$ by
\[\oct(\gamma)\{r-s\}(P) \defeq \sum_{e \in v\rightarrow \gamma v}\kappa_{\f}\{r-s\}(e,P),\]
where $v \rightarrow \gamma v$ denotes the (unique) path from $v$ to $\gamma v$ in the tree and where $\kappa_{\f}$ is the modular symbol on the tree from Section \ref{modular symbols on the tree}. From the properties of $\kappa_{\f}$ in Proposition \ref{modular symbol tree properties}, we see that this is a cocycle in $\mathrm{Z}^1(\Gamma,\Delta_{k,k})$ and its class $\oc \in \HGammas$ is independent of $v$.
\end{mdef}

We have natural actions of the Hecke operators on $\HGammas$, as in \cite{HidE}, Chapter 6.3. As in \cite{Ort04}, using Proposition \ref{properties of double integral} and Lemma \ref{modular symbol tree properties}, we have the following.
\begin{mprop}
Both $\oc$ and $\lc$ are eigensymbols for the Hecke operators away from $\n$ with the same eigenvalues as $f$, that is, they are elements in the $f$-isotypic subspace of $\HGammas$.
\end{mprop}
\end{shortversion}

\begin{longversion}

\subsection{Definition of $\lc$}
For ease of notation, let $\Delta_{k,k} \defeq \mathrm{Hom}(\Delta_0,V_{k,k}(\Cp)^*).$
\begin{mdef}\label{define lc}
Let $\tau \in \uhp_{\pri}$, and define a function $\Gamma \rightarrow \Delta_{k,k}$ by
\[\lct(\gamma)\{r-s\}(P) \defeq \int_{\tau}^{\gamma \tau}\int_r^s (P).\]
\end{mdef}
From the properties of the double integral above, we see that this is a cocycle in $\mathrm{Z}^1(\Gamma,\Delta_{k,k})$ and its class $\lc \in \HGammas$ is independent of $\tau$.

\subsection{Definition of $\oc$}
\begin{mdef}Let $v \in \vertices(\treep)$, and define a function $\Gamma \rightarrow \Delta_{k,k}$ by
\[\oct(\gamma)\{r-s\}(P) \defeq \sum_{e \in v\rightarrow \gamma v}\kappa_{\f}\{r-s\}(e,P),\]
where $v \rightarrow \gamma v$ denotes the (unique) path from $v$ to $\gamma v$ in the tree and where $\kappa_{\f}$ is the modular symbol on the tree from Section \ref{modular symbols on the tree}.
\end{mdef}
From the properties of $\kappa_{\f}$ in Proposition \ref{modular symbol tree properties}, we see that this is a cocycle in $\mathrm{Z}^1(\Gamma,\Delta_{k,k})$ and its class $\oc \in \HGammas$ is independent of $v$.

\subsection{Hecke operators}
We have natural actions of the Hecke operators on $\HGammas$, as in \cite{HidE}, Chapter 6.3. As in \cite{Ort04}, using Proposition \ref{properties of double integral} and Lemma \ref{modular symbol tree properties}, we have the following.
\begin{mprop}
For $I$ coprime to $\n$, we have
\[\lc|T_I = \lambda_I \lc\]
and 
\[\oc|T_I = \lambda_I \oc.\]
In particular, both $\lc$ and $\oc$ are elements of the $f$-isotypic subspace of $\HGammas$.
\end{mprop}

\end{longversion}

%
%

\section{Relating to $L$-values}

\subsection{Optimal embeddings}
In both \cite{Dar01} and \cite{Ort04}, one obtains a link between the cohomology classes defined above and $L$-values of the form $f$ at characters of conductor $c$ by exploiting \emph{optimal embeddings of conductor $c$}, that is, certain embeddings $F\times F\hookrightarrow \mathrm{M}_2(F)$. Heuristically, one should see such an embedding as giving rise to a cycle in the symmetric space corresponding to our level group, and hence it gives rise to a homology class that we can then pair our cohomology classes against (see, for example, \cite{GMS15} for details of this approach). 
$\lb$
In practice we can work out the details of this rather more explicitly. In particular, each such embedding gives rise to an infinite cyclic subgroup $\langle \gamma\rangle$ of $\Gamma$, which has precisely two fixed points $x,y$ in $\Proj(F_{\pri})$. One then writes down a polynomial $P(t)$ that is invariant under the action of $\gamma$, and evaluates cohomology classes in $\HGamma$ at $\gamma, x, y$ and $P$. In this section, we bypass embeddings completely, and simply write down the corresponding values of $\gamma, x, y$ and $P$. Throughout, and for the remainder of the paper, we will assume $k$ is even (in lieu of Corollary \ref{exceptional zero existence}).
\begin{mdef}
Let $c \in \roi_F$ be prime to $\pri = (\pi)$ and $v\in \roi_F$ prime to $c$.
\begin{itemize}
\item[(i)] Define $s'\defeq \text{order of }\pi\text{ in }(\roi_F/c)^\times,$ and let $s \defeq 2\times(\text{order of }\pi^2\text{ in }(\roi_F/c)^\times).$
\item[(ii)] Define
\[\gamma_{c,v} \defeq \matrd{\pi^{-s/2}}{(\pi^{s/2}-\pi^{-s/2})\frac{v}{c}}{0}{\pi^{s/2}}\]
(As $c|(\pi^{s/2}-\pi^{-s/2})$, there is no problem with denominators). Note that the fixed points of this in $\Proj(F_{\pri})$ are $-v/c$ and $\infty$.
\item[(iii)] Let $M_{c,v}(t) \defeq ct + v$, a M\"{o}bius transformation taking $-v/c$ to $0$ and fixing $\infty$.
\item[(iv)] Let $P_{c,v}(t) \defeq (ct+v)^{k/2}(\overline{ct}+\overline{v})^{k/2},$ an element of $V_{k,k}(L)$ that is invariant under the action of $\gamma_{c,v}$.
\end{itemize}
\end{mdef}
\begin{longversion}
\begin{mrem}
Again, the choice of action on modular symbols means we have to be careful about how $\gamma_{c,v}$ is acting on the various spaces involved here. We have an action of $\gamma_{c,v}$ on $\Proj(F)$ by the usual fractional linear transformations $x \mapsto (ax+b)/(cx+d)$, and under this action, $v/c$ and $\infty$ are fixed. We also have an action of $\gamma_{c,v}$ on $\uhp_{\pri}$ by $\tau \mapsto (b+d\tau)/(a+c\tau)$, and under this action, for any choice of $\tau \in \uhp_{\pri}$, we have
\[\lim_{j \rightarrow \infty} \gamma_{c,v}^j\tau = -v/c, \hspace{12pt} \lim_{j\rightarrow \infty}\gamma_{c,v}^{-j}\tau = \infty.\]
\end{mrem}
\end{longversion}
\begin{shortversion}
\begin{mrem}
The actions of $\gamma_{c,v}$ on the cusps and on $\Proj(F_{\pri})$ and $\uhp_{\pri}$ are different. In particular $\gamma_{c,v}$ fixes the cusps $v/c$ and $\infty$, whilst for any $\tau \in \uhp_{\pri}$, we have $\lim_{j \rightarrow \infty} \gamma_{c,v}^j\tau = -v/c$ and $\lim_{j\rightarrow \infty}\gamma_{c,v}^{-j}\tau = \infty.$
\end{mrem}
\end{shortversion}
\begin{longversion}
The following is a collection of well-known and/or easily proved facts and properties that we'll need when relating cohomology classes to $L$-values.
\end{longversion}
\begin{mprop}\label{embedding invariants}
\begin{itemize}
\item[(i)] Let $t\in\Proj(F_{\pri})$. Then we have $M_{c,v}(\gamma_{c,v} t) = \pi^s M_{c,v}(t)$.
\item[(ii)] There is a bijection between doubly infinite (non-repeating) paths in $\treep$ and distinct pairs of elements of $\Proj(F_{\pri})$. (Here by \emph{doubly infinite} we mean that every vertex in the path is connected to precisely two other vertices, that is, the path has no `end').
\item[(iii)] Let $\mathrm{path}(-v/c,\infty)$ denote the infinite path in $\treep$ corresponding to the pair $(-v/c,\infty)$. It is possible to index the vertices and edges of this path such that $t(e_j) = v_j$ and 
\[U(e_j) = \{t \in \Proj(F_{\pri}): v_{\pri}(M_{c,v}(t)) \geq -j\}.\]
\item[(iv)] Under the action of $\Gamma$ on $\treep$, we have $\gamma_{c,v}e_j = e_{j+s}.$
\item[(v)] Let $U(v_j) \defeq \{t \in \Proj(F_{\pri}): v_{\pri}(M_{c,v}(t)) = - j\}.$ Then  
\[F_{c,v} \defeq \bigcup_{j=0}^{s-1}U(v_{j})\]
is a fundamental domain for the action of $\gamma_{c,v}$ on $\Proj(F_{\pri})\backslash\{-v/c,\infty\}$.

\end{itemize}
\end{mprop}
\begin{proof}
Part (i) is direct calculation. For part (ii), see \cite{DT08}. Parts (iii) to (v) follow from the same arguments as the corresponding statements in \cite{Dar01} over $\Q$, which are given in equations (88) for part (iii), (90) for part (iv) and (87) for part (v).
\end{proof}

\subsection{Relating $\mathrm{oc}_f$ to classical $L$-values}\label{ocf lvalues}
With the above given, we are in a position to compute a special value of the cocycle $\oct$. 
\begin{mprop}
Let $J_v$ denote the coset $v\langle\pi\rangle \subset (\roi_F/c)^\times.$ For each $a \in J_v$, let $j(a)$ be the smallest non-negative integer such that $a \equiv v\pi^{j(a)} \newmod{c}$. We have
\[\oct(\gamma_{c,v})\{v/c - \infty\}(P_{c,v}) = \beta\sum_{a\in J_v} \omega^{j(a)}\phi_f\{a/c - \infty\}\bigg[(ct+a)^{k/2}(\overline{ct}+\overline{a})^{k/2}\bigg],\]
where
\[\beta = \left\{\begin{array}{ll} 1 &: s = s',\\
2 &: \omega = 1, s = 2s',\\
0 &: \omega = -1, s = 2s'.
\end{array}\right.\]
\end{mprop}
\begin{proof}
We compute directly that
\[\oct(\gamma_{c,v})\{v/c-\infty\}(P_{c,v}) = \sum_{j=1}^s\kappa_{\f}\{x-y\}(e_j,P_{c,v}).\]
Let $u\in\roi_F$ be any element such that $u \equiv v/c \newmod{p^s}$. Since $U(e_j) = \{t\in \Proj(F_{\pri}): \pi^{j}(ct+v) \in \roi_{\pri}\} = \{t \in \Proj(F_{\pri}): \pi^{j}(t+u)\in\roi_{\pri}\} = U(\gamma_je_*),$ where 
\begin{equation}\label{gammaj}
\gamma_j \defeq \matrd{\pi^{-j}}{u}{0}{1} \in \Gammat,
\end{equation}
we see that $e_j = \gamma_j e_*$. Hence $\kappa_f\{v/c-\infty\}(e_j,P_{c,v}) = \omega^j\phi_f|\gamma_j^{-1}\{v/c-\infty\}(P_{c,v}).$ Now direct computation shows that 
\[\gamma_j^{-1}(v/c) = \frac{\pi^{j}v - cu\pi^{j}}{c},\hspace{12pt} \gamma_j^{-1}\infty = \infty,\]
and
\[\gamma_j^{-1}\cdot P_{c,v}(t) = (ct + (\pi^{j}v - cu\pi^{j}))^{k/2}(\overline{ct} + (\overline{\pi}^{j}\overline{v} - \overline{cu}\overline{\pi}^{j}))^{k/2}. \]
As $j$ ranges from 1 to $s$, the quantity $a = \pi^{j}v-cu\pi^{j}$ ranges over $J_v$, once if $s = s'$ and twice if $s = 2s'$. In the latter case, the signs $\omega^j$ at each instance are both equal to 1 if $\omega = 1$ and are opposite if $\omega = -1$. Finally, the result follows from the observation that if $a,a' \in \roi_F$ such that $a \equiv a' \newmod{c}$, then the action of $\smallmatrd{1}{(a-a')/c}{0}{1} \in \Gamma_0(\n)$ shows that 
\[\phi_f\{a/c-\infty\}[(ct+a)^{k/2}(\overline{ct}+\overline{a})^{k/2})] = \phi_f\{a'/c-\infty\}[(ct+a')^{k/2}(\overline{ct}+\overline{a'})^{k/2})].\qedhere\]
\end{proof}
Henceforth we assume we are not in the case $\beta = 0$, which will be irrelevant to our results. An immediate corollary is the following.
\begin{mcor}\label{oc and L-values}
Suppose $\chi$ is a finite order character of conductor $c$ (that we naturally see as a character on $(\roi_F/c)^\times$), and that $\chi(\pri) = \omega$. Then 
\[\sum_{v \in (\roi_F/c)^\times}\chi(v)\oc(\gamma_{c,v})\{v/c - \infty\}(P_{c,v}) = s\left[\frac{D^{r+1}\tau(\chi^{-1})|c|^r \unitsize}{2\Omega_f}\right] \Lambda(f,\chi,k/2+1).\]
\end{mcor} 
\begin{proof}
The corresponding result for the cocycle $\oct$ follows from equation (\ref{integralformula}). In particular, note that the condition on $\chi(\pi)$ means that $\chi(v)\omega^{j(a)} = \chi(v\pi^{j(a)}) = \chi(a)$, whilst the sum runs over all cosets $v\langle\pi\rangle \subset (\roi_F/c)^\times$ for all $v\in(\roi_F/c)^\times$, which means that every element of $(\roi_F/c)^\times$ is hit $s'$ times. But $\beta s' = s$.
$\lb$
To see that the result holds at the level of cohomology classes, let $b \in \mathrm{B}^1(\Gamma,\Delta_{k,k})$ be a coboundary. Then there exists $\phi \in \Delta_{k,k} =\mathrm{Hom}(\Delta_0,V_{k,k}(\Cp)^*)$ such that $b(\gamma)\{r-s\}(P) = (\phi|\gamma)\{r-s\}(P) - \phi\{r-s\}(P).$ But it follows directly that $b(\gamma_{c,v})\{v/c-\infty\}(P_{c,v}) = 0$ using the $\gamma_{c,v}$-invariance of $v/c, \infty$ and $P_{c,v}$.
\end{proof}

\subsection{Relating $\mathrm{lc}_f$ to $p$-adic derivatives}
In this section, we prove a formula analogous to the one above, this time relating special values of the class $\lc$ to the derivative of the $p$-adic $L$-function at the corresponding critical value. For clarity, we will assume throughout this section that $p$ is split; the cases where $p$ is inert or ramified are almost identical (and, in fact, slightly easier where they differ). We freely use the notation of Section \ref{ocf lvalues}.

\subsubsection{Rephrasing as an integral over the fundamental domain}
Let $c$ and $v$ be as before. The following lemma will be useful in the sequel.

\begin{mlem} \label{integral zero}We have 
\[\int_{F_{c,v}\times \mathcal{O}_{\overline{\mathfrak{p}}}}P_{c,v}(t)d\mu_{\mathcal{F}}\left\lbrace \frac{v}{c}- \infty \right\rbrace (t)=0.\]
\end{mlem}
\begin{longversion}
\begin{proof} Note that from Proposition \ref{embedding invariants} we have $F_{c,v}= U(e_{s-1})\smallsetminus U(e_{-1})$ and $\gamma_{c,v}e_{-1}= e_{s-1}$. Hence $U(e_{s-1})= \gamma_{c,v}^{-1}U(e_{-1})$, so that $U(e_{s-1})\times \mathcal{O}_{\overline{\mathfrak{p}}}= \gamma_{c,v}^{-1}\left[ U(e_{-1})\times \mathcal{O}_{\overline{\mathfrak{p}}}\right]$. From Proposition 5.11 we deduce that 
\begin{align*}
\int_{F_{c,v}\times \mathcal{O}_{\overline{\mathfrak{p}}}}P_{c,v}(t)d\mu_{\mathcal{F}}&\left\lbrace \frac{v}{c}- \infty \right\rbrace (t) = 
	\int_{U(e_{s-1})\times \mathcal{O}_{\overline{\mathfrak{p}}}}P_{c,v}(t)d\mu_{\mathcal{F}}\left\lbrace \frac{v}{c}- \infty \right\rbrace (t) - \int_{U(e_{-1})\times \mathcal{O}_{\overline{\mathfrak{p}}}}P_{c,v}(t)d\mu_{\mathcal{F}}\left\lbrace \frac{v}{c}- \infty \right\rbrace (t)\\
&=	 \int_{\gamma_{c,v}^{-1}\left[U(e_{-1})\times \mathcal{O}_{\overline{\mathfrak{p}}}\right]}P_{c,v}(t)d\mu_{\mathcal{F}}\left\lbrace \frac{v}{c}- \infty \right\rbrace (t) - \int_{U(e_{-1})\times \mathcal{O}_{\overline{\mathfrak{p}}}}P_{c,v}(t)d\mu_{\mathcal{F}}\left\lbrace \frac{v}{c}- \infty \right\rbrace (t) \\
	&= \int_{U(e_{-1})\times \mathcal{O}_{\overline{\mathfrak{p}}}}\gamma_{c,v}^{-1}\cdot P_{c,v}(t)d\mu_{\mathcal{F}}\left\lbrace \gamma_{c,v}^{-1}\left(\frac{v}{c}\right) - \infty \right\rbrace (t) - \int_{U(e_{-1})\times \mathcal{O}_{\overline{\mathfrak{p}}}}P_{c,v}(t)d\mu_{\mathcal{F}}\left\lbrace \frac{v}{c}- \infty \right\rbrace (t)\\
&= 0,
\end{align*}
where	 the last equality follows from the fact that $P_{c,v}(t)$ and $\frac{v}{c} \in \mathbb{P}^{1}(F)$ are invariant under the action of $\gamma_{c,v}$.
\end{proof}
\end{longversion}
\begin{shortversion}
\begin{proof}
This is a direct calculation using the fact that $F_{c,v} = U(e_{s-1})\backslash U(e_{-1})$ and $U(e_{s-1}) = \gamma_{c,v}^{-1}U(e_{-1})$, combined with the invariance of $P_{c,v}$, $c/v$ and $\infty$ under $\gamma_{c,v}$.
\end{proof}
\end{shortversion}

\begin{longversion}
Recall that by definition,
\[\tilde{\mathrm{lc}}_{f, \tau}(\gamma_{c,v})\{v/c - \infty \}(P_{c,v})= \int_{\mathbb{P}^{1}_{\mathfrak{p}}}\mathrm{log}_{p}\left( \frac{t- \gamma_{c,v}\tau}{t- \tau}\right)P_{c,v}(t)d\mu_{\mathcal{F}}\left\lbrace \frac{v}{c}- \infty \right\rbrace (t).\]
To ease notation, we'll henceforth write $\mu_{\f}^{c,v}$ for $\mu_{\f}\{v/c-\infty\}$. In the following proposition, we simplify the definition.
\end{longversion}
\begin{shortversion}
Recall the definition of $\widetilde{\mathrm{lc}}_{f,\tau}(\gamma_{c,v})\{v/c - \infty\}(P_{c,v})$ from Definition \ref{define lc}. We rephrase this in a way that allows us to see $p$-adic $L$-values. To ease notation, we'll henceforth write $d\mu_{\f}^{c,v}$ for $d\mu_{\f}\{v/c-\infty\}$. 
\end{shortversion}
\begin{mprop}\label{integral over fund domain}
We have
$$\lct(\gamma_{c,v})\{v/c - \infty \}(P_{c,v})= \int_{F_{c,v}\times \mathcal{O}_{\overline{\mathfrak{p}}}}\mathrm{log}_{p}(ct+ v)P_{c,v}(t)d\mu_{\mathcal{F}}^{c,v}(t).$$
\begin{longversion}In particular, the left hand side does not depend on the choice of $\tau$.\end{longversion}
\end{mprop}
The proof will require a lemma. We have a decomposition
\[\Proj_{\pri} = \left(U(\overline{e}_{(n+1)s-1})\sqcup U(e_{-ns-1}) \bigsqcup_{i = -n}^n\gamma_{c,v}^i F_{c,v}  \right)\times \roi_{\pribar}.\]
We break the integral up into sums over these components and let $n$ tend to $\infty$.
\begin{mlem}\label{endpoints zero}
The integral over the `endpoints' vanishes in the limit. More precisely, let
\[I(n) \defeq \int_{U'(\overline{e}_{(n+1)s-1})\sqcup U'(e_{-ns-1})}\mathrm{log}_{p}\left( \frac{t- \gamma_{c,v}\tau}{t- \tau}\right)P_{c,v}(t)d\mu_{\mathcal{F}}^{c,v}(t).\]
Then $\lim_{n\rightarrow \infty}I(n) = 0.$ In particular, 
\[\tilde{\mathrm{lc}}_{f, \tau}(\gamma_{c,v})\{v/c - \infty \}(P_{c,v}) =\lim_{n\rightarrow \infty}\sum_{i= -n}^{n}\int_{\gamma_{c,v}^{i}\left[ F_{c,v}\times \mathcal{O}_{\overline{\mathfrak{p}}}\right] }\mathrm{log}_{p}\left( \frac{t- \gamma_{c,v}\tau}{t- \tau}\right)P_{c,v}(t)d\mu_{\mathcal{F}}^{c,v}(t).\]
\end{mlem}
\begin{proof}
Since $U'(\overline{e}_{(n+1)s-1}) = \gamma_{c,v}^{-n-1}U'(e_{-1})$ and $U'(e_{-ns-1}) = \gamma_{c,v}^nU'(e_{-1}),$ we can write
\begin{align*}
\int_{U'(\overline{e}_{(n+1)s-1})}\mathrm{log}_{p}&\left(\frac{t- \gamma_{c,v}\tau}{t- \tau}\right)P_{c,v}(t)d\mu_{\mathcal{F}}^{c,v}(t) =\\
& -\int_{U'(e_{-1})}\mathrm{log}_{p}\left( \frac{\gamma_{c,v}^{n+1}t- \gamma_{c,v}\tau}{\gamma_{c,v}^{n+1}t- \tau} \right)P_{c,v}(t)d\mu_{\mathcal{F}}^{c,v}(t)
\end{align*}
using the invariance of $P_{c,v}$ and $d\mu_{\f}^{c,v}$ (and similarly for the integral over $U'(e_{-ns-1})$). For each fixed $n$, the expressions
\[\frac{\gamma_{c,v}^{n+1}t- \gamma_{c,v}\tau}{\gamma_{c,v}^{n+1}t- \tau}, \hspace{12pt} \frac{t- \gamma_{c,v}^{-n}\tau}{t- \gamma_{c,v}^{-n-1}\tau}\]
are both M\"{o}bius transformations sending $\gamma_{c,v}^{-n}\tau$ to zero and $\gamma_{c,v}^{-n-1}\tau$ to $\infty$, hence one is a constant scalar multiple of the other. A short calculation using Lemma \ref{integral zero} now shows that
\[I(n) = \int_{U(e_{-1})}\mathrm{log}_p\left(\frac{t- \gamma_{c,v}^{-n}\tau}{t- \gamma_{c,v}^{-n-1}\tau}\cdot \frac{t- \gamma_{c,v}^{n+1}\tau}{t- \gamma_{c,v}^{n}\tau}\right) P_{c,v}(t)d\mu_{\f}^{c,v}(t).\]
In the limit, the expression in the $\mathrm{log}_p$ becomes $M_{c,v}(\gamma_{c,v}t)/M_{c,v}(t).$ But by Proposition \ref{embedding invariants}, this is equal to $\pi^s$, and $\mathrm{log}_p(\pi^s) = 0$. Hence the integrand tends to zero and the first result follows. The second result is a direct consequence of this and the above decomposition.
\end{proof}

\begin{shortversion}
\begin{proof}\emph{(Proposition \ref{integral over fund domain})}.

We use the second part of Lemma \ref{endpoints zero}. A general term of the sum is an integral over $\gamma_{c,v}^i[F_{c,v}\times\roi_{\pribar}],$ and using the transformation property, we can write it as an integral over $F_{c,v}\times\roi_{\pribar}$. In particular, using similar methods to the proof of Lemma \ref{endpoints zero}, we see that
\[\lct(\gamma_{c,v})\left\{\frac{v}{c} - \infty\right\}(P_{c,v}) = \lim_{n\rightarrow\infty} \sum_{i=-n}^n \int_{F_{c,v}\times\roi_{\pribar}} \mathrm{log}_p\left(\frac{t-\gamma_{c,v}^{i+1}\tau}{t-\gamma_{c,v}^i\tau}\right)P_{c,v}(t)d\mu_{\f}^{c,v}(t).\]
This sum telescopes via the expression inside the $\mathrm{log}_p$, and in the limit, the resulting expression tends to $M_{c,v}(t) = ct+v$. This completes the proof.
\end{proof}

\end{shortversion}

\subsubsection{Rephrasing using the overconvergent modular symbol}
We can decompose the fundamental domain of Proposition \ref{integral over fund domain} into a union of sets of the form $U(v_{j})\times\roi_{\pribar}$, and hence as 
\[F_{c,v}\times\roi_{\pribar} = \bigsqcup_{j=0}^{s-1}\gamma_j^{-1}\left[\roi_{\pri}^\times\times\roi_{\pribar}\right],\]
where $\gamma_j$ is as defined in equation (\ref{gammaj}). Here we've used the fact that $U(v_j) = U(\gamma_j v_*) = \gamma_j^{-1}\roi_{\pri}^\times,$ which can be seen directly. Hence we have
\begin{align}\label{lcf value}
\lct(\gamma_{c,v})&\{v/c-\infty\}(P_{c,v}) = \sum_{j=0}^{s-1}\int_{\gamma_j^{-1}[\roi_{\pri}^\times\times\roi_{\pribar}]} \mathrm{log}_p(ct+v)P_{c,v}(t)d\mu_{\f}^{c,v}\\
&= \sum_{j=0}^{s-1}\omega^j\int_{\roi_{\pri}^\times\times\roi_{\pribar}}\mathrm{log}_p(ct+a_j)(ct+a_j)^{k/2}(\overline{ct}+\overline{a_j})^{k/2}d\mu_{\f}\{a_j/c-\infty\},\notag
\end{align}
where $a_j = \pi^j v - cu\pi^j$ (using the calculations and notation of the previous section). \begin{shortversion}Now a short explicit calculation shows:\end{shortversion}
\begin{mlem}\label{lcf oms}
Write $c = g\overline{\pi}^r$, for some $r\geq 0$ and $(g)$ coprime to $(p)$. For any $a\in\roi_F$, we have 
\begin{align*}\int_{\roi_{\pri}^\times\times\roi_{\pribar}}\mathrm{log}_p(ct+a)&(ct+a)^{k/2}(\overline{ct}+\overline{a})^{k/2}d\mu_{\f}^{c,v}(t) \\
&= (g\overline{g})^{k/2}\lambda_{\pribar}^{r}\Psi_f\left|\matrd{1}{a}{0}{g}\right.\{0-\infty\}\bigg(\mathrm{log}_p(t)(t\overline{t})^{k/2}\ind_{a\newmod{\pribar^r}}\bigg)\\
&= \lambda_{\pribar}^r\int_{\{[a]\}\times\roi_\pri^\times\times\big[a+\pribar^r\roi_{\pribar}\big]}\mathrm{log}_p(z_{\pri})\langle\mathbf{z}_p\rangle^{k/2}d\mu_p(\mathbf{z}),
\end{align*}
where $\Psi_f$ is the overconvergent modular symbol attached to $f$, $\ind_{a\newmod{\pri^r}}$ is the indicator function of the set $\roi_{\pri}^\times \times[a + \pribar^r\roi_{\pribar}],$ and $\mu_p$ is the distribution on $(\roi_F/\mathfrak{g})^\times\times(\roi_F\otimes_{\Z}\Zp)$ defined in Section \ref{padic lfunction section}.
\end{mlem}

\begin{longversion}
\begin{proof}
Since $\roi_{\pri}^\times\times\roi_{\pribar} \subset \roi_F\otimes_{\Z}\Zp$, we can write the integral as 
\[\Psi_f\{a/c-\infty\}(\mathrm{log}_p(ct+a)(ct+a)^{k/2}(\overline{ct}+\overline{a})^{k/2}\ind_{\roi_{\pri}^\times\times\roi_{\pribar}}).\]
In the right hand side, replace $\Psi_f$ with $\lambda_{\pribar}^{-r}\Psi_f|U_{\pribar}^r$ and expand out the sum define the Hecke operator. Since $U_{\pri}$ commutes with the action of $\smallmatrd{1}{a}{0}{g}$, this gives
\begin{align*}
\mathrm{(RHS)} &= (g\pi^r\overline{g\pi}^r)^{k/2}\sum_{b \newmod{\pribar^r}}\Psi_f\left|\matrd{1}{a}{0}{g}\matrd{1}{b}{0}{\overline{\pi}^r}\right.\{0-\infty\}\bigg(\mathrm{log}_p(t)(t\overline{t})^{k/2}\ind_{a\newmod{\pribar^r}}\bigg)\\
&= (c\overline{c})^{k/2}\Psi_f\left|\matrd{1}{a}{0}{c}\right.\{0-\infty\}\bigg(\mathrm{log}_p(t)(t\overline{t})^{k/2}\ind_{a\newmod{\pribar^r}}\bigg),
\end{align*}
since the indicator kills all other terms of the sum. Expanding this out gives the left hand side. The second equality is by definition of $\mu_p$.
\end{proof}
\end{longversion}

\begin{mcor}\label{lcf oms cor}We have
\[\lc(\gamma_{c,v})\{v/c-\infty\}(P_{c,v}) = \sum_{a \in J_v}\omega^{j(a)}\lambda_{\pribar}^r\int_{\{[a]\}\times\roi_\pri^\times\times\big[a+\pribar^r\roi_{\pribar}\big]}\mathrm{log}_p(z_{\pri})\langle\mathbf{z}_p\rangle^{k/2}d\mu_p(\mathbf{z}).\]
\end{mcor}
\begin{proof}
This follows directly from equation (\ref{lcf value}), Lemma \ref{lcf oms} and the fact that as $j$ ranges from $0$ to $s-1$, $a_j$ ranges over the set $J_v$ defined in the previous section with multiplicity $\beta$. As in the previous section, this value depends only on the cohomology \emph{class}, so we drop the tilde.
\end{proof}

\subsubsection{Completing the calculation}
Finally, we tie the above results together to relate $\mathrm{lc}_{f}$ to derivatives of the $p$-adic $L$-function attached to $f$. As in the last subsection we suppose that we are not in the case $\beta= 0$.

\begin{mcor}\label{lc and L-values}
Suppose $\chi$ is a finite order character of conductor $c$ (that we naturally see as a character on $(\roi_F/c)^\times$), and that $\chi(\pri) = \omega$. Write $c = g\overline{\pi}^r$. Then 
 \[\left[1 - \frac{\chi(\pribar)N(\pribar)^{k/2}}{\lambda_{\pribar}}\right]\sum_{v \in \left(\roi_{F}/ c\right)^{\times}}\chi(v)\lc(\gamma_{c,v})\{v/c - \infty \}(P_{c,v})= s\lambda_{\pribar}^r \frac{d}{ds_{\mathfrak{p}}}L_{p}(f, \chi, \mathbf{s})\mid_{\mathbf{s}= \mathbf{k/2}}.\]
\end{mcor}
\begin{proof} Note that if $\mathbf{z} \in \cl_F(\mathfrak{g}p^\infty)$ satisfies $\mathbf{z}\equiv \pi^{j}v \newmod{c}$ then $\chi(\mathbf{z})= \omega^{j}\chi(v)$. We use Corollary \ref{lcf oms cor} and the fact that $s = \beta s'$. In particular,  we see that
\begin{align*}\sum_{v \in (\roi_F/ c)^{\times}}&\chi(v)\lc(\gamma_{c,v})\{v/c - \infty \}(P_{c,v})
= \\
&s\lambda_{\pribar}^r\int_{(\roi_F/\mathfrak{g})^\times \times[\roi_{\pri}^\times \times \roi_{\pribar}]}\chi(\mathbf{z})\mathrm{log}_p(z_{\pri})\langle\mathbf{z}_p\rangle^{k/2}d\mu_p(\mathbf{z}).
\end{align*}
We need an integral over the units in $\roi_{\pribar}$ in order to descend to $\cl_F(\mathfrak{g}p^\infty)$. To this end, note that we can break the integral up as 
\[\int_{(\roi_F/\mathfrak{g})^\times \times[\roi_{\pri}^\times \times \roi_{\pribar}]} = \int_{(\roi_F/\mathfrak{g})^\times \times[\roi_{\pri}^\times \times \roi_{\pribar}^\times]} + \int_{(\roi_F/\mathfrak{g})^\times \times[\roi_{\pri}^\times \times \pribar\roi_{\pribar}]}.\]
If $\pribar|c$, then extend $\chi$ to $\roi_{\pribar}$ by defining it to be zero on $\pribar\roi_{\pribar}$ (and correspondingly on the larger space in the integral). Then
\begin{align*}\int_{(\roi_F/\mathfrak{g})^\times \times[\roi_{\pri}^\times \times \pribar\roi_{\pribar}]}&\chi(\mathbf{z})\mathrm{log}_p(z_{\pri})\langle \mathbf{z}_p\rangle^{k/2}d\mu_p(\mathbf{z}) \\
&= \lambda_{\pribar}^{-1}\int_{(\roi_F/\mathfrak{g})^\times \times[\roi_{\pri}^\times \times \roi_{\pribar}]}\chi(\overline{\pi}\mathbf{z})\mathrm{log}_p(\overline{\pi}z_{\pri})\langle \overline{\pi}\mathbf{z}_p\rangle^{k/2}d\mu_p(\mathbf{z})\\
&= \frac{\chi(\overline{\pi})N(\pribar)^{k/2}}{\lambda_{\pribar}}\int_{(\roi_F/\mathfrak{g})^\times \times[\roi_{\pri}^\times \times \roi_{\pribar}]}\chi(\mathbf{z})\mathrm{log}_p(z_{\pri})\langle \mathbf{z}_p\rangle^{k/2}d\mu_p(\mathbf{z}).
\end{align*}
By collecting these terms together, and using the fact that the integrand is invariant under units, we see that
\begin{align*}\left[1 - \frac{\chi(\pribar)N(\pribar)^{k/2}}{\lambda_{\pribar}}\right]\sum_{v \in (\roi_F/ c)^{\times}}\chi(v)&\lc(\gamma_{c,v})\{v/c - \infty \}(P_{c,v})\\
&= s\lambda_{\pribar}^r\int_{\cl_F(\mathfrak{g}p^\infty)}\chi(\mathbf{z})\mathrm{log}_p(z_{\pri})\langle\mathbf{z}_p\rangle^{k/2}d\mu_p(\mathbf{z}).
\end{align*}
But the latter integral is precisely the derivative in the $s_{\pri}$ direction of the $p$-adic $L$-function (as can be seen by differentiating the integrand in the definition). This completes the proof.
\end{proof}

%
%
\section{The cohomology of $\Gamma$}
We now prove that the $f$-isotypic component of $\HGamma$, that is the subspace of classes $\psi \in \HGamma$ such that $T_I \psi = \lambda_I \psi$ for all $I$ coprime to $\n$, is one dimensional. We make use of the tree.
\begin{mprop}\label{long exact sequence}
Let $\Delta_{k,k} \defeq \HomD$. There is an exact sequence
\begin{align*}\big[\h^0(\Gamma_0(\m),\Delta_{k,k})\oplus\h^0(\Gamma'_0(\m),\Delta_{k,k}) \big]\longrightarrow &\hspace{2pt}\h^0(\Gamma_0(\n),\Delta_{k,k}) \longrightarrow \h^1(\Gamma,\Delta_{k,k}) \\
&\longrightarrow \big[\h^1(\Gamma_0(\m),\Delta_{k,k})\oplus\h^1(\Gamma'_0(\m),\Delta_{k,k}))\big]\end{align*}

where $\Gamma'_0(\m) \defeq \alpha\Gamma_0(\m)\alpha^{-1}.$
\end{mprop}
\begin{proof}
This is \cite{Ser80}, Chapter II.2.8, Proposition 13, applied to 
\[M = \HomD\]
and $G = \Gamma,$ noting that $\mathrm{Stab}_{\Gamma}(e_*) = \Gamma_0(\n)$, $\mathrm{Stab}_{\Gamma}(s(e_*)) = \Gamma_0(\m)$, and $\mathrm{Stab}_{\Gamma}(t(e_*)) = \Gamma'_0(\m)$.
\end{proof}

\subsection{A cohomological lemma}
We can relate these spaces to compactly supported cohomology using a trick of Ash and Stevens (see \cite{AS86}, Proposition 4.2).
\begin{mlem}\label{ash stevens}
 For each $i\geq 0$ there is a Hecke-equivariant isomorphism
	$$\mathrm{H}^{i}(\Gamma_{0}(\mathfrak{m}), \Delta_{k,k})) \cong \mathrm{H}^{i+1}_{c}(\Gamma_{0}(\mathfrak{m})\backslash \mathcal{H}_{3}, \widetilde{V_{k,k}(\C_{p})^{\ast}}),$$
where the tilde denotes the natural local system attached to the right $\Gamma_{0}(\mathfrak{m})$-module $V_{k,k}(\C_{p})^{\ast}$. We have identical isomorphisms for $\Gamma_{0}'(\mathfrak{m})$ and $\Gamma_{0}(\mathfrak{n})$.
\end{mlem}

\begin{longversion}
In fact, this is a special case of a more general result.  We briefly recall the construction of the Borel--Serre compactification of the manifold $\Gamma_0(\mathfrak{m})\backslash \mathcal{H}_{3}$. In \cite{BS74} Borel and Serre construct a space $\overline{\mathcal{H}}_{3}$, containing $\mathcal{H}_{3}$, and prove that it is a manifold with corners and smooth boundary. More precisely we have 
\[\overline{\mathcal{H}}_{3}= \mathcal{H}_{3} \sqcup \bigsqcup_{x \in \Proj(F)}e_{x},\]
where for each $x \in \Proj(F)$ we define $e_{x}\cong \mathbb{R}^{2}$. Moreover, $\overline{\mathcal{H}}_{3}$ is a contractible space. The action of $\Gamma_0(\mathfrak{m})$ on $\mathcal{H}_{3}$ extends to an action on $\overline{\mathcal{H}}_{3}$ and the quotient $\Gamma_0(\mathfrak{m})\backslash \overline{\mathcal{H}}_{3}$ is in fact a compact manifold with corners and it is called the \emph{Borel--Serre compactification} of $\Gamma_0(\mathfrak{m})\backslash \mathcal{H}_{3}$. 
\end{longversion}
\begin{shortversion}
In fact, this is a special case of a more general result. Recall first the following definition.
\begin{mdef}
The \emph{Borel--Serre compactification} of $\Gamma_0(\mathfrak{m})\backslash\uhs$ is $\Gamma_0(\mathfrak{m})\backslash\overline{\uhp}_3,$ where 
\[\overline{\uhp}_3 = \uhs \sqcup \bigsqcup_{x\in\Proj(F)}e_x,\] 
with each $e_x$ a copy of $\R^2$ (see the appendix of \cite{Ser70} or \cite{BS74} for details).
\end{mdef}
\end{shortversion}
Let $M$ be a $\GLt(F)$-module.\footnote{More generally, we may replace $\GLt(F)$ with a semigroup $\Lambda\subset\GLt(F)$ containing $\Gamma_0(\m)$ and the matrices defining the Hecke operators.} Restricting the action to $\Gamma_0(\m)$, we obtain natural local systems on the space $\Gamma_0(\m)\backslash\uhs$, its compactification $\Gamma_0(\m)\backslash\overline{\uhs}$, and the boundary $\partial(\Gamma_0(\mathfrak{m})\backslash \overline{\mathcal{H}}_{3})$; in an abuse of notation, we write $\widetilde{M}$ for each of these local systems. For any $i\geq 0$ we have the following \begin{shortversion}Hecke equivariant\end{shortversion} isomorphisms:  
 \begin{equation} \label{group vs manifold 1}
 \rmH^{i}\bigg(\Gamma_0(\mathfrak{m}), \rmH^{0}\big(\overline{\mathcal{H}}_{3}, M\big)\bigg) \cong \rmH^{i}\left(\Gamma_0(\mathfrak{m}) \backslash \mathcal{H}_{3}, \widetilde{M}\right);
 \end{equation}
 \begin{equation}\label{group vs manifold 2}
 \rmH^{i}\bigg(\Gamma_0(\mathfrak{m}), \rmH^{0}\big(\partial\overline{\mathcal{H}}_{3}, M\big)\bigg) \cong \rmH^{i}\left(\partial\big[\Gamma_0(\mathfrak{m}) \backslash \overline{\mathcal{H}}_{3}\big], \widetilde{M}\right);
 \end{equation}
 \begin{equation}\label{group vs manifold 3}
 \rmH^{i}\bigg(\Gamma_0(\mathfrak{m}), \rmH^{1}\big(\overline{\mathcal{H}}_{3}, \partial \overline{\mathcal{H}}_{3}, M\big)\bigg) \cong \rmH_{c}^{i+1}\left(\Gamma_0(\mathfrak{m}) \backslash \mathcal{H}_{3}, \widetilde{M}\right).
 \end{equation}
\begin{longversion}
Using the action of $\GLt(F)$ on $M$ we can define Hecke operators on the cohomology groups on the left hand side of each isomorphism above, and we can prove that these isomorphisms are in fact Hecke equivariant.
\end{longversion}
\begin{mlem}
For each $i\geq 0$, the following diagram commutes:
\[\xymatrix@C=0.7em@R=1.7em{ \ar[r]& \rmH^{i}(\Gamma_0(\mathfrak{m}), \mathrm{Hom}(\Delta_{0}, M)) \ar[r]\ar^{\wr}[d]&  \rmH^{i+1}(\Gamma_0(\mathfrak{m}),M)  \ar[r]\ar^{\wr}[d]&  \rmH^{i+1}(\Gamma_0(\mathfrak{m}), \mathrm{Hom}(\Delta, M)) \ar[r]\ar^{\wr}[d]& \\
	\ar[r]&  \rmH^{i+1}_{c}(\Gamma_0(\mathfrak{m})\backslash \mathcal{H}_{3}, \widetilde{M}) \ar[r]& \rmH^{i+1}(\Gamma_0(\mathfrak{m})\backslash \mathcal{H}_{3}, \widetilde{M}) 	\ar[r]&  \rmH^{i+1}(\partial\left(\Gamma_0(\mathfrak{m})\backslash \overline{\mathcal{H}}_{3}\right) , \widetilde{M}) \ar[r]& .}\]
Here the vertical maps are isomorphisms, the horizontal sequences are exact, and each map is Hecke equivariant. Moreover, we have analogous diagrams for $\Gamma_{0}'(\mathfrak{m})$ and $\Gamma_{0}(\mathfrak{n})$.
\end{mlem}	
\begin{proof} The proof follows the same strategy as \cite{AS86}.	Let $\Delta\defeq \mathrm{Div}(\Proj(F))$. We know that the variety $\overline{\mathcal{H}}_{3}$ is contractible and its boundary, $\partial\overline{\mathcal{H}}_{3}$, is a disjoint union of copies of $\mathbb{R}^{2}$ indexed by $\Proj(F)$. Hence we have the following isomorphisms:
\[\rmH_{0}(\overline{\mathcal{H}}_{3}, \Z)\cong \Z \ \ \ \mathrm{and} \ \ \ \rmH_{0}(\partial\overline{\mathcal{H}}_{3}, \Z)\cong \Delta.\]
Additionally, the boundary map 
$\rmH_{1}(\overline{\mathcal{H}}_{3}, \partial\overline{\mathcal{H}}_{3}, \Z) \rightarrow \rmH_{0}(\partial\overline{\mathcal{H}}_{3}, \Z)$ 
induces an isomorphism 
$\rmH_{1}(\overline{\mathcal{H}}_{3}, \partial\overline{\mathcal{H}}_{3}, \Z)\cong\Delta_{0}.$ 
Moreover, all these isomorphisms are equivariant with respect to the action of $\mathrm{GL}_{2}(F)$, and assemble into a commutative diagram
\[
\xymatrix@R=1.7em{ 
	0 \ar[r]& \rmH_{1}(\overline{\mathcal{H}}_{3}, \partial\overline{\mathcal{H}}_{3}\, \Z) \ar[r]\ar^{\wr}[d]& \rmH_{0}(\partial\overline{\mathcal{H}}_{3}, \Z)\ar[r]\ar^{\wr}[d]& \rmH_{0}(\overline{\mathcal{H}}_{3}, \Z) \ar[r]\ar^{\wr}[d]& 0\\
	0 \ar[r]&  \Delta_0  \ar[r]&  \Delta  \ar[r]&  \Z  \ar[r]& 0.
}
\]
Taking $\mathrm{Hom}(\bullet, M)$ in the diagram above we obtain the following commutative diagram where the horizontal sequences are exact and the vertical arrows are isomorphisms:
\begin{equation}
\label{before cohomology}
\xymatrix@R=1.7em{ 
	0 \ar[r]& \rmH^{0}(\overline{\mathcal{H}}_{3}, M) \ar[r]\ar^{\wr}[d]& \rmH^{0}(\partial\left(\overline{\mathcal{H}}_{3}\right), M)  \ar[r]\ar^{\wr}[d]& \rmH^{1}(\overline{\mathcal{H}}_{3}, \partial\overline{\mathcal{H}}_{3}, M)  \ar[r]\ar^{\wr}[d]& 0\\
	0 \ar[r]&  M  \ar[r]&  \mathrm{Hom}(\Delta, M)  \ar[r]&  \mathrm{Hom}(\Delta_0 , M)  \ar[r]& 0.
}
\end{equation}
The diagram above is compatible with the action of $\GLt(F)$ on each term. Considering the long exact sequence of the cohomology of $\Gamma_0(\mathfrak{m})$ attached to each line in the diagram above and using isomorphisms (\ref{group vs manifold 1}), (\ref{group vs manifold 2}) and (\ref{group vs manifold 3}) completes the proof. Since each diagram respects the action of $\GLt(F)$, the diagram of the lemma is Hecke equivariant.   
\end{proof}		
\begin{longversion}
\begin{mrem}Note that this also gives an isomorphism 
\[\symb_{\Gamma_0(\n)}(V_{k,k}(L)^*) \cong \h^1_{\mathrm{c}}\big(\Gamma_0(\n)\backslash\uhs, \widetilde{V_{k,k}(L)^*}\big).\]
This cohomological setting gives the appropriate way of generalising modular symbols to other settings.
\end{mrem}
\end{longversion}

\subsection{Multiplicity one for group cohomology}
\begin{mthm}[Multiplicity one]\label{multiplicity one}
Let $f \in S_{k,k}(\Omega_0(\n))^{\pri-\mathrm{new}}$. Then the system of Hecke eigenvalues corresponding to $f$ appears in the cohomology groups 
\[\h^i_{\mathrm{c}}(\Gamma_0(\n)\backslash\uhs,\widetilde{V_{k,k}(\Cp)^*})\]
(only) for $i = 1,2$, and the $f$-isotypic component of these spaces is one-dimensional. The system of eigenvalues corresponding to $f$ does not appear in the cohomology of $\Gamma_0(\m)\backslash\uhs$ or $\Gamma'_0(\m)\backslash\uhs$.
\end{mthm}
\begin{proof}
For the first statment, see \cite{Hid94}, Proposition 3.1, which is given for cuspidal cohomology. The cuspidal cohomology can be viewed as a subspace of the compactly supported cohomology, and cokernel of the inclusion is an Eisenstein subspace, that is, the $f$-isotypic parts of $\h^i_{\mathrm{cusp}}$ and $\h^i_{\mathrm{c}}$ are the same (see \cite{Har87}, Section 3.2.5, or \cite{Urb95}, Theorem 3.2). The second statement holds simply because $f$ is new at $\pri$, so its system of eigenvalues cannot appear at level $\m$.
\end{proof}
The main result of this section is the following:
\begin{mthm}\label{cohomology dimension}
Let $H$ be any space with an action of the Hecke operators $T_I$ for $I$ coprime to $\n$, and let $H_{(f)}$ denote the $f$-isotypic component.
\begin{itemize}
\item[(i)] We have $\h^i(\Gamma_0(\m),\Delta_{k,k})_{(f)} = \h^i(\Gamma'_0(\m),\Delta_{k,k})_{(f)} = 0$, for $i = 0,1$.
\item[(ii)] The space $\h^0(\Gamma_0(\n),\Delta_{k,k})_{(f)}$ is a one-dimensional $\Cp$-vector space.
\item[(iii)] The space $\h^1(\Gamma, \Delta_{k,k})_{(f)}$ is a one-dimensional $\Cp$-vector space.
\end{itemize}
\end{mthm}
\begin{proof}
Parts (i) and (ii) are immediate from Lemma \ref{ash stevens} and Theorem \ref{multiplicity one}. A direct calculation shows that the maps in the long exact sequence of Proposition \ref{long exact sequence} are equivariant with respect to the Hecke operators, and hence we get a corresponding exact sequence on the $f$-isotypic parts. But from parts (i) and (ii), the end terms are zero whilst the second term is one-dimensional. Part (iii) follows.
\end{proof}
\begin{mrem}
An identical proof would show that the $f$-isotypic component of $\h^2(\Gamma,\Delta_{k,k})$ is also one-dimensional. 
\end{mrem}

%
%
\section{Exceptional zeros of Bianchi modular forms}
In this section, we draw our previous results together to prove statements about exceptional zeros.
\begin{mlem}
The cohomology class $\oc$ is non-zero.
\end{mlem}
\begin{proof}We exploit a result of Rohrlich on the non-vanishing of $L$-values, namely the theorem in the introduction of \cite{Roh91}. If $\oc$ vanishes, then by Corollary \ref{oc and L-values} so must the $L$-values $L(f,\chi,k/2+1)$ for all Dirichlet characters $\chi$ with $\chi(\pi) = \omega$, and this contradicts Rohrlich's result.
\end{proof}
\begin{mdef}
Define $\mathcal{L}_{\pri} \in \Cp$ to be the unique scalar such that
\[\lc = \mathcal{L}_{\pri}\oc \in \h^1(\Gamma,\Delta_{k,k})_{(f)},\]
which is possible by Theorem \ref{cohomology dimension}.
\end{mdef}
The following is the main result of this paper.
\begin{mthm}
Let $f \in S_{k,k}(\Omega_0(\n))$ be a Bianchi eigenform over an imaginary quadratic field of class number 1, and suppose that $f$ is new at $\pri = (\pi)$ and has small slope at every other prime above $p$. Let $-\omega$ be the eigenvalue of the Atkin--Lehner operator $W_{\pri}$, and let $\chi$ be a finite order Hecke character with $\chi(\pri) = \omega$. Let $\ff$ be the $p$-part of the conductor of $\chi$. Then there is an exceptional zero $L_p(f,\chi,\mathbf{k/2}) = 0$, and there exists $\mathcal{L}_{\pri} \in \Cp$ such that
\begin{align*}
\frac{d}{ds_{\pri}}&L_p(f,\chi\tm^{k/2},\mathbf{s})\big|_{\mathbf{s} = \mathbf{k/2}} =\\
& \mathcal{L}_{\pri}\left(\prod_{\substack{\mathfrak{q}|p\\ \mathfrak{q}\neq \pri}}Z_{\mathfrak{q}}(\chi,k/2)\right)\left[\frac{D^{k/2+1}\tau(\chi^{-1})|c|^r \unitsize}{\Omega_f\lambda_{\ff}}\right] \Lambda(f,\chi,k/2+1).
\end{align*}
Moreover, $\mathcal{L}_{\pri}$ is independent of the choice of $\chi$ satisfying the above conditions.
\end{mthm}
\begin{proof}
The existence of an exceptional zero was Corollary \ref{exceptional zero existence}. The result follows from Corollaries \ref{oc and L-values} and \ref{lc and L-values} and the definition of $\mathcal{L}_{\pri}$.
\end{proof}

%
%
\section{An arithmetic description in the base-change case}
Let $f$ be a cuspidal Bianchi modular form of weight $(k,k)$. In the case where $f$ is the base-change of a classical cusp form $\tilde{f}$ (of weight $k+2$), new at $p$, then we are able to relate the $\mathcal{L}$-invariants at $p$ of $f$ and $\tilde{f}$ under a non-vanishing assumption. Indeed, Artin formalism says that there is a relation
\[\Lambda(f,s) = \Lambda(\tilde{f},s)\Lambda(\tilde{f},\chi_{F/\Q},s),\]
where $\chi_{F/\Q}$ is the quadratic character associated to the extension $F/\Q$. 
\begin{longversion}We can extend this to be true of the values of $p$-adic $L$-functions at critical values.
\begin{mlem}
Let $\chi$ be a (rational) Dirichlet character of conductor $cp^n$, where $c$ is coprime to $p$ and $n\geq 1$. Then we have
\[\left[\frac{D^{2(r+1)}(cp^n)^{2(r+1)}}{\tau(\chi\cdot N_{F/\Q})\tau(\chi_F)\lambda_{p^n}(f)}\right]\Lambda(f,\chi,r+1) = \left[\frac{(cp^n)^{r+1}}{\tau(\chi)\lambda_{p^n}(\tilde{f})}\right]\Lambda(\tilde{f},\chi,r+1)\left[\frac{(cp^nD)^{r+1}}{\tau(\chi\chi_F)\lambda_{p^n}(\tilde{f})}\right]\Lambda(\tilde{f},\chi\chi_F,r+1).\]
In particular, under suitable normalisations, we have
\[(2/ \unitsize)L_p(f,(\chi\cdot N_{F/\Q})\tm^r,\mathbf{r}) = L_p(\tilde{f},\chi,\mathbf{r}\tm^r, \mathbf{r})L_p(\tilde{f},\chi\chi_F\tm^r,\mathbf{r}).\]
\end{mlem} 
\begin{proof}
Most of the constant terms are evidently equal, and the equality of the $L$-function terms is Artin formalism. The only things left to check are the eigenvalues and the Gauss sums. But $\lambda_p(f) = \lambda_p(\tilde{f})^2$ by inspection, whilst we have $\tau(\chi\cdot N_{F/\Q}) = \tau(\chi)^2\tau(\chi_F)$ by an explicit calculation combined with a characteristic zero version of the classical Hasse--Davenport identity for Gauss sums. If we let $\Omega_{\tilde{f}}^\pm$ be the periods of $\tilde{f}$, then $\Omega_f \defeq \Omega_{\tilde{f}}^+\Omega_{\tilde{f}}^-\tau(\chi_F)$ is a period of $f$, and using these periods, the statement about $p$-adic $L$-values follows immediately from the respective interpolation formulae.
\end{proof}
Given an equality at critical values, we can lift this result to $p$-adic $L$-functions: 
\end{longversion}
\begin{shortversion}We can lift this to the level of $p$-adic $L$-functions.
\end{shortversion}
\begin{mprop}
Suppose $f$ is the base-change to $F$ of a classical modular form that has small slope at $p$. Given periods $\Omega_{\tilde{f}}^\pm$ for $\tilde{f}$, we may choose the period of $f$ to be $\Omega_f = \Omega_{\tilde{f}}^+\Omega_{\tilde{f}}^-\tau(\chi_{F/\Q})$. The natural inclusion $\Zp \hookrightarrow \roi_F\otimes_{\Z}\Zp$ allows us to restrict the resulting $p$-adic $L$-function of $f$ to an analytic function on $\Zp$ (the \emph{cyclotomic variable}). This restriction factorises as
\[L_p(f,s) = \frac{ \unitsize}{2}\cdot L_p(\tilde{f},s)\cdot L_p(\tilde{f},\chi_{F/\Q},s)\]
for all $s \in \Zp$.
\end{mprop}
\begin{proof}
\begin{shortversion}
The $p$-adic $L$-functions of $f$ and $\tilde{f}$ are defined using canonical overconvergent modular symbols $\Psi_f$ and $\Psi_{\tilde{f}}$ respectively, which give rise to locally analytic distributions $\mu_p$ on $\cl_F(p^\infty)$ and $\tilde{\mu}_p$ on $\Zp^\times$ respectively (see \cite{PS11} for the classical case). We can restrict $\mu_p$ to functions that factor through the norm map and hence consider it as a distribution on $\Zp^\times.$ Finally, define a locally analytic distribution $\tilde{\mu}^{\chi_{F/\Q}}_p$ on $\Zp^\times$ as follows: use $\Psi_{\tilde{f}}$ to define a distribution $\tilde{\mu}_{p,D'}$ on $(\Z/D')^\times \times \Zp^\times$ as in Section \ref{padic lfunction section}, where $D'$ is the prime to $p$ part of $D$, and define
\[\tilde{\mu}_p^{\chi_{F/\Q}}(\zeta) = \tilde{\mu}_{p,D'}\left(\chi_{F/\Q}\cdot\zeta\right),\]
where we consider $\chi_{F/\Q}\cdot\zeta$ as a function on $(\Z/D')^\times\times\Zp^\times$ in the natural way.
\end{shortversion}
\begin{longversion}
Recall the definition of $L_p(f,s)$ from Definition \ref{padiclfunctiondefinition}: we defined a distribution $\mu_p$ on $\cl_F(p^\infty)$ and constructed an analytic function on $\roi_F\otimes_{\Z}\Zp$ using it. The $p$-adic $L$-function of $\tilde{f}$ can be constructed similarly (see \cite{PS11}, for example) using a distribution $\tilde{\mu}_p$ on $\cl_{\Q}^+(p^\infty) \cong \Zp^\times$. We can restrict $\mu_p$ to the characters that factor through the norm map $N_{F/\Q} : \cl_F(p^\infty) \rightarrow \cl_{\Q}^+(p^\infty)$ to obtain a distribution on $\Zp^\times$. Finally, define a locally analytic distribution $\tilde{\mu}_p^{\chi_{F/\Q}}$ on $\Zp^\times$ as follows: extend $\tilde{\mu}_p$ to a distribution $\tilde{\mu}_{p,D'}$ on $(\Z/D')^\times \times\Zp^\times$ as in Section \ref{padic lfunction section}, where $D'$ is the prime to $p$ part of $D$, and define 
\[\tilde{\mu}_p^{\chi_{F/\Q}}(\zeta) = \tilde{\mu}_{p,D'}\left(\chi_{F/\Q}\cdot\zeta\right),\]
where we consider $\chi_{F/\Q}\cdot\zeta$ as a function on $(\Z/D')^\times\times\Zp^\times$ in the natural way.
\end{longversion}
$\lb$
The result follows if we can prove the (stronger) result that for any locally analytic function $\zeta$ on $\Zp^\times$, we have
\[\mu_p(\zeta\circ N_{F/\Q}) = \frac{ \unitsize}{2}\cdot\tilde{\mu}_p(\zeta)\cdot\tilde{\mu}_p^{\chi_{F/\Q}}(\zeta).\]
For $\zeta$ of the form $\chi|\cdot|^r$, for $0 \leq r \leq k$ and $\mathrm{cond}(\chi)|p^\infty$, it is a lengthy but elementary\footnote{The most difficult part of the calculation is to show that the respective Gauss sum terms agree when $p$ is inert; this is a characteristic zero analogue of the Hasse--Davenport identity.} check that this holds by the respective interpolation formulae and the relation between the classical $L$-functions of $f$ and $\tilde{f}$ above. But both of these distributions are admissible of order $k < k+1$, and hence equality on this set of points forces equality everywhere, and the proof is complete.
\end{proof}

Now suppose that $f$ is new (hence small slope) at $p$, so $L_p(f,s)$ has an exceptional zero at $s = k/2$, and that $p$ is inert in $F$. In this case, $\chi_{F/\Q}(p) = -1$, and accordingly we see that
$L_p(\tilde{f},s)$ has an exceptional zero at $s = k/2$ whilst $L_p(\tilde{f},\chi_{F/\Q},s)$ does not. Differentiating, we see that
\[\frac{d}{ds}L_p(f,s)\big|_{s=k/2} = \frac{ \unitsize}{2}\left[\frac{d}{ds}L_p(\tilde{f},s)\big|_{s=k/2}\right]L_p(\tilde{f},\chi_{F/\Q},k/2).\]
Using the respective interpolation conditions (for the $p$-adic $L$-functions of both $f$ and $\tilde{f}$), 
\begin{longversion}collapsing the exceptional factor in the second classical term, and using classical Artin formalism, we find that 
\begin{align*}\mathcal{L}_{\pri}\cdot\left[\frac{|D|^{k/2+1}w}{\Omega_f}\right] \Lambda(f,k/2+1) = \mathcal{L}_p(\tilde{f})\cdot\left[\frac{2|D|^{k/2+1}}{\Omega_{\tilde{f}}^+\Omega_{\tilde{f}}^-\tau(\chi_{F/\Q})}\right]\cdot\Lambda(f,k/2+1),
\end{align*}
where $\mathcal{L}_p(\tilde{f})$ denotes the $\mathcal{L}$-invariant of $\tilde{f}$ at $p$. Without loss of generality, we may renormalise so that $\Omega_f = \Omega_{\tilde{f}}^+\Omega_{\tilde{f}}^-\tau(\chi_{F/\Q})$. Then we deduce the following:
\end{longversion}
\begin{shortversion}
we deduce the following:
\end{shortversion}
\begin{mprop}
Suppose $f$ is the base-change of a classical modular form $\tilde{f}$ to $F$, and that $p$ is inert in $F$. If $\Lambda(f,k/2+1) \neq 0$, we have
\[\mathcal{L}_{\pri} = 2\mathcal{L}_p(\tilde{f}),\]
where $\mathcal{L}_p(\tilde{f})$ is the $\mathcal{L}$-invariant of $\tilde{f}$.
\end{mprop}
A similar argument shows that if $p$ ramifies in $F$, then $\mathcal{L}_{\pri} = \mathcal{L}_p(\tilde{f})$ (the difference arising when we consider the exceptional factor in the second classical term, which equals 2 in the inert case and 1 in the ramified case). More generally, since the $\mathcal{L}$-invariant is purely local, similar methods show that these results also hold if there exists a Dirichlet character $\chi$ over $F$, factoring through the norm to $\Q$, such that $\Lambda(f,\chi,k/2+1) \neq 0$.
\begin{mrems}
\begin{itemize}
\item[(i)]In \cite{Tri06}, Trifkovi\'{c} uses local methods to construct conjectural global points on modular elliptic curves over imaginary quadratic fields. He associates a (conjectural) weight $(0,0)$ Bianchi modular form $f$ to such a curve and defines the points using double integrals of the sort considered in this paper. His main conjecture is Conjecture 6, but this is only well-defined given his previous Conjecture 5, which is the statement that a value of the double integral generates a lattice commensurable with $q^{\Z}$ (where $q$ is the Tate period of the elliptic curve). In the analogous case over $\Q$, this is known via the results of \cite{GS93}, but it remains open over imaginary quadratic fields in general. When $f$ is base-change and $p$ is inert or ramified, the above proposition combined with \cite{GS93} gives a proof of Trifkovi\'{c}'s Conjecture 5.
\item[(ii)] Results of a similar nature for abelian base-change in the case of certain weight 2 automorphic forms for $\GLt$ were recently obtained by Gehrmann in \cite{Geh17}.
\end{itemize}
\end{mrems}
%
%
\section{Further remarks}
It is natural to expect exceptional zero results of the form above to exist in wider generality. We conclude by commenting on possible extensions of the results of this paper.

\subsection{Higher class number}\label{class number}
The most obvious restriction we've assumed in the work above is that the class number of $F$ is 1. This is not a serious restriction to the method, but rather a choice made for ease of exposition; indeed, the authors carried out a large portion of the work without assuming any restriction on the class number at all, but found that the increase in technicalities made the arguments much longer and harder to follow. To generalise whilst retaining classical methods, one would need to consider functions on $h$ copies of the tree, where $h$ is the class number, and find elements in the direct sum of $h$ group cohomology groups, and the Hecke operators permute these functions/classes.

\subsection{Higher derivatives}
In the case where $p$ is split, it is possible that even $(d/ds_{\pri})L_p(f,\chi,\mathbf{s})|_{\mathbf{s} = \mathbf{k/2}} = 0$ even when the corresponding classical $L$-value is non-zero, since we are left with an exceptional factor at the other prime $\pribar$ that could vanish. The most natural solution to this would be to differentiate again, this time in the $s_{\pribar}$-direction. In considering this approach, the authors (directly) constructed classes $\mathrm{oc}_{f},\mathrm{lc}_f \in \h^2(\Gamma,\Delta_{k,k})$ that again both live in a one-dimensional eigenspace, and hoped to relate their special values to $L$-values. However, whilst it seemed likely that such results would exist at the level of \emph{cocyles}, there was a fundamental obstacle to this approach in that these special values were no longer well-defined on the corresponding cohomology classes. Results of this nature do exist, however, when the weight of $f$ is 2 (parallel weight 0); over totally real fields, this was done by Spie{\ss} (see \cite{Spi14}), and more recently over arbitrary number fields by Deppe and Bergunde (see \cite{Dep16} and \cite{Ber17}). Their approach, which whilst cohomological is quite different to that considered here, is to construct higher cohomology classes by taking cup products of cycles in $\h^1$. It would be interesting to generalise their results to arbitrary weights.

\subsection{Arithmetic descriptions of $\mathcal{L}_{\pri}$ in general}
We have already seen that it is possible to obtain an arithmetic description of the $\mathcal{L}$-invariant in the case that $f$ is the base-change of a classical modular form. Evidently it is desirable to obtain a similar description when $f$ is not base-change. This case, however, appears far harder. The difficulty arises from the absence of nice geometry, which is a fundamental block to generalising classical methods of determining the $\mathcal{L}$-invariant in arithmetic terms. One of the most interesting cases to study is when the form is weight $(0,0)$ (also referred to as weight 2) and corresponds, via the (conjectural) imaginary quadratic version of modularity, to an elliptic curve over $F$. In this situation, the elliptic curve will have multiplicative reduction at $\pri$, so there is a Tate uniformiser $q \in F_{\pri}^\times$ such that $E(\overline{F}_{\pri}) \cong \overline{F}_{\pri}^\times/q^{\Z}.$ By analogy with the case of Hilbert modular forms, one would expect that:
\begin{mcon}\label{exp desc}Let $Q_{\pri} \defeq \mathrm{Norm}_{F_{\pri}/\Qp}(q),$
and let $f_{\pri} \defeq [\roi_{\pri}/\pri:\Fp]$ denote the inertia degree of $\pri$.
Then we have
\[\mathcal{L}_{\pri} = f_{\pri}\frac{\mathrm{log}_{p}(Q_{\pri})}{\mathrm{ord}_{p}(Q_{\pri})}.\]
\end{mcon}
(Compare \cite{Mok09}, Theorem 8.2 and \cite{Hid09}, Theorem 1.2). The computational results of Trifkovi\'{c} in\cite{Tri06}, and later of Guitart, Masdeu and \c{S}eng\"{u}n in \cite{GMS15}, can be seen as providing strong evidence for this conjecture. More generally, for a form $f$ of arbitrary weight, let $\rho_{f}$ be the Galois representation attached to $f$ by Taylor (see \cite{Tay94}). A natural generalisation of Conjecture \ref{exp desc} is:
\begin{mcon} The $\mathcal{L}$-invariant $\mathcal{L}_{\mathfrak{p}}$ is the Benois--Greenberg $\mathcal{L}$-invariant attached to $\rho_{f}$ (see \cite{Ben11} and \cite{Ros15}). 
\end{mcon}
The authors hope to address these questions in future work.

\small
\renewcommand{\refname}{\normalsize References} 
\bibliography{references}{}
\bibliographystyle{alpha}
\Addresses

\end{document}